\documentclass{article}

\usepackage[margin=1.25in]{geometry}
\usepackage{hyperref}
\usepackage{natbib}
\usepackage{graphicx} 
\usepackage{subfigure} 
\usepackage{longtable}
\usepackage{bm}
\usepackage{algorithm}
\usepackage{algorithmic}
\usepackage{booktabs}
\usepackage{amsmath,amssymb,amsfonts}
\usepackage{multirow}
\newtheorem{theorem}{Theorem}

\newtheorem{lemma}{Lemma}
\newtheorem{corollary}{Corollary}

\newcommand{\x}{\mathbf{x}}
\newcommand{\uu}{\mathbf{u}}
\newcommand{\y}{\mathbf{y}}
\newcommand{\z}{\mathbf{z}}
\newcommand{\w}{\mathbf{w}}
\renewcommand{\b}{\mathbf{b}}
\newcommand{\X}{\mathbf{X}}

\newcommand{\A}{\mathbf{A}}

\newcommand{\I}{\mathbf{I}}

\newcommand{\G}{\mathbf{G}}

\newcommand{\E}{\mathbb{E}}
\newcommand{\hlam}{\hat{\bm\lambda}}
\newcommand{\tlam}{\tilde{\bm\lambda}}

\newcommand{\olam}{\bm\lambda}
\newcommand{\tx}{\mathbf{\tilde{x}}}
\newcommand{\tna}{\tilde{\nabla}}
\newcommand{\tb}{\tilde{\b}}

\DeclareMathOperator*{\argmin}{argmin}

\begin{document}

\title{Faster and Non-ergodic $O(1/K)$ \\Stochastic  Alternating Direction Method of Multipliers}

\author{\\
Cong Fang\\
\texttt{fangcong@pku.edu.cn}\\
Peking University \\
\and \\
Feng Cheng\\
\texttt{fengcheng@pku.edu.cn}\\
Peking University \\
\and \\
Zhouchen Lin\\
\texttt{zlin@pku.edu.cn} \\
Peking University \\
}

%

%
%
%
%
%
\date{Original circulated date: 22th April, 2017.}
\maketitle

\begin{abstract}
We study stochastic convex optimization subjected to linear equality constraints. Traditional Stochastic Alternating Direction Method of Multipliers~\cite{STOC-ADMM} and its Nesterov's acceleration scheme~\cite{OPT-SADMM} can only achieve  ergodic $O(1/\sqrt{K})$ convergence rates, where $K$ is the number of iteration. By introducing Variance Reduction~(VR) techniques, the convergence rates improve to ergodic   $O(1/K)$~\cite{SAG-ADMM,SVRG-ADMM}. In this paper, we propose a new stochastic ADMM  which elaborately integrates   Nesterov's extrapolation and VR techniques. We prove that  our algorithm can achieve a non-ergodic $O(1/K)$ convergence rate which is optimal for separable linearly constrained non-smooth convex problems, while  the convergence rates of VR based ADMM methods are  actually tight  $O(1/\sqrt{K})$ in non-ergodic sense.  To the best of our knowledge, this is the \emph{first} work that achieves a truly accelerated, stochastic convergence rate for constrained convex problems.  The experimental results demonstrate that our  algorithm is faster than the existing state-of-the-art stochastic ADMM methods.
\end{abstract} 
\section{Introduction}
We consider the following general convex finite-sum problems with linear constraints:
\begin{eqnarray}\label{problem}
\min_{\x_1,\x_2} &&h_1(\x_1)+f_1(\x_1) +h_2(\x_2)+\frac{1}{n}\sum_{i=1}^n f_{2,i}(\x_2),\notag\\
s.t. && \A_1 \x_1 +\A_2 \x_2 = \b,
\end{eqnarray}
where  $f_{1}(\x_1)$ and  $f_{2,i}(\x_2)$ with $i\in\{1,2,\cdots,n\}$  are convex and have Lipschitz continuous gradients,  $h_1(\x_1)$ and $h_2(\x_2)$ are also  convex.  We denote that $L_1$ is the Lipschitz constant of $f_{1}(\x_1)$,   $L_2$ is the Lipschitz constant of $f_{2,i}(\x_2)$ with $i\in\{1,2,\cdots,n\}$, and $f_2(\x)=\frac{1}{n}\sum_{i=1}^n f_{2,i}(\x)$. We define $F_i(\x_i) =h_i(\x_i)+ f_i(\x_i)$ for $i = 1,2$, $\x=(\x_1,\x_2)$, $F(\x_1,\x_2)=F_1(\x_1)+F_2(\x_2)$, and $\A\x = \sum_{i=1}^2\A_i\x_i$.

Problem~\eqref{problem} is of great importance in machine learning. The finite-sum function $f_2(\x_2)$ is typically a loss over training samples, and the remaining functions control the structure or regularize the model to aid generalization~\cite{OPT-SADMM}. The idea of using linear constraints to decouple the loss and  regularization terms enables researchers to consider some more sophisticated regularization terms which might be very complicated to solve through proximity operators for Gradient Descent~\cite{FISTA} methods. For example, for multitask learning problems~\cite{multitask,shenli}, the regularization term is set as $\mu_1\|\x\|_* +\mu_2 \|\x\|_1$, for most graph-guided fused Lasso and overlapping group Lasso problem ~\cite{lasso,SVRG-ADMM}, the regularization term can be written as $\mu  \|\A \x \|_1$,  and for many multi-view learning tasks~\cite{Wang2016Joint},  the regularization terms always involve $\mu_1\|\x\|_{2,1} +\mu_2 \|\x\|_* $. 

Alternating Direction Method of Multipliers~(ADMM) is a very popular optimization method to solve Problem~\eqref{problem}, with its advantages in speed, easy implementation, good scalability, shown in lots of literatures~(see survey \cite{boyd2011distributed}). However, though ADMM is effective in practice, the provable convergence rate is not fast. A popular criterion to judge convergence is in ergodic sense. And it is proved in ~\citep{he20121,LADMPSAP} that ADMM converges with an $O(1/K)$ ergodic rate. Since the non-ergodic results~($\x^K$), rather than the ergodic one~(convex combination of $\x^1,\x^2,\cdots, \x^K$) is much faster in practice, researchers gradually  turn to analyse the convergence rate in non-ergodic sense. \citet{ADMM} prove that the Douglas-Rachford~(DR) splitting converges in  non-ergodic $O(1/\sqrt{K})$. They also construct a family of functions showing that non-ergodic $O(1/\sqrt{K})$ is tight.  ~\citet{chen2015inertial} establish $O(1/\sqrt{K})$ for Linearized ADMM. Then \citet{LADM-NE} accelerate ADMM through Nesterov's extrapolation and obtain a non-ergodic $O(1/K)$ convergence rate. They also prove that the lower complexity bound of ADMM type methods for the separable linearly constrained nonsmooth convex problems is exactly $O(1/K)$, which demonstrates that their algorithm is optimal. The convergence rates for different ADMM based algorithms are shown in Table \ref{convergence rate}.

\begin{table}[tb]
	\centering
	\label{convergence rate}
	\caption{Convergence rates of ADMM type methods solving Problem~\eqref{problem} ~(``non-'' indicates ``non-ergodic'', while ``er-'' indicates ``ergodic''. ``Sto.'' is short for ``Stochstic'',  and ``Bat.'' indicates batch or deterministic algorithms).}
	\begin{tabular} {|c |c |c|}\hline
		\multirow{1}*{Type} &Algorithm & Convergence Rate \\\hline\hline
		\multirow{2}*{\centering Bat.}& ADMM~\scriptsize{\citep{ADMM}}        &  Tight non-$O(\frac{1}{\sqrt{K}})$\\\cline{2-3}
		&LADM-NE~\scriptsize{\citep{LADM-NE}}         & Optimal non-$O(\frac{1}{K})$ \\\hline\hline
		\multirow{7}*{\centering Sto.}&STOC-ADMM~\scriptsize{\citep{STOC-ADMM}}          &  er-$O(\frac{1}{\sqrt{K}})$ \\\cline{2-3}
		&OPG-ADMM~\scriptsize{\citep{OPG-ADMM}}         & er-$O(\frac{1}{\sqrt{K}})$ \\\cline{2-3}
		&OPT-ADMM~\scriptsize{\citep{OPT-SADMM}}         & er--$O(\frac{1}{\sqrt{K}})$\\\cline{2-3}
		&SDCA-ADMM~\scriptsize{\citep{SDCA-ADMM}}         & unknown \\\cline{2-3}
		&SAG-ADMM~\scriptsize{\citep{SAG-ADMM}}         & Tight non--$O(\frac{1}{\sqrt{K}})$ \\\cline{2-3}
		&SVRG-ADMM~\scriptsize{\citep{SVRG-ADMM}}         & Tight non--$O(\frac{1}{\sqrt{K}})$ \\\cline{2-3}
		&ACC-SADMM~\scriptsize{(ours)}         & Optimal non--$O(\frac{1}{K})$ \\\hline
	\end{tabular}
	\label{resultss}
\end{table}

On the other hand, to meet the demands of solving  large-scale machine learning problems, stochastic algorithms~\cite{bottou2004stochastic} have drawn a lot of interest in recent years. For stochastic ADMM~(SADMM), the prior works are from STOC-ADMM~\cite{STOC-ADMM} and OPG-ADMM~\cite{OPG-ADMM}. Due to the noise of gradient, both of the two algorithms can only achieve an ergodic $O(1/\sqrt{K})$ convergence rate.  There are two lines of research to accelerate SADMM. The first is to introduce the Variance Reduction~(VR)~\cite{SVRG,SAGA,SAG} techniques into SADMM. VR methods are widely accepted tricks  for finite sum problems which ensure the descent direction to have a bounded variance and so can achieve faster convergence rates.  The existing VR based SADMM algorithms include SDCA-ADMM~\cite{SDCA-ADMM}, SAG-ADMM~\cite{SAG-ADMM} and  SVRG-ADMM~\cite{SVRG-ADMM}. SAG-ADMM and SVRG-ADMM can provably  achieve   ergodic $O(1/K)$ rates for Porblem~\eqref{problem}. However,  in  non-ergodic  sense,  their convergence  rates are $O(1/\sqrt{K})$~(please see detailed discussions in Section~\ref{VR-SADMM}), while the fastest rate for batch ADMM method is $O(1/K)$~\cite{LADM-NE}. So there still exists a gap between stochastic and batch~(deterministic) ADMM. The second line to accelerate SADMM is through the Nesterov's acceleration~\cite{nesterov1983method}. This work is from \cite{OPT-SADMM}, in which the authors  propose an ergodic $O(\frac{R^2}{K^2}+\frac{D_y+\rho}{K}+\frac{\sigma}{\sqrt K})$ stochastic algorithm~(OPT-ADMM).  The dependence on the smoothness constant of the convergence rate is $O(1/K^2)$ and  so each term in the convergence rate seems to have been improved to optimal. This method is imperfect due to the following two reasons: 1) The worst convergence  rate  is still the ergodic $O(1/ \sqrt K)$. There is no theoretical improvement in the order $K$. 2)  OPT-SADMM  is  not very effective in practice.  The method does not adopt any special technique to tackle the noise of  gradient except adding a proximal term $\frac{\|\x^{k+1}-\x^{k}\|^2}{\sigma k^{3/2}}$  to ensure convergence. As the gradients have noise,  directly applying the original Nesterov's extrapolation on the variables  often causes the  objective function  to oscillate and decrease slowly during iteration.

In this paper, we propose Accelerated Stochastic ADMM (ACC-SADMM) for large scale general convex finite-sum problems with linear constraints.  By elaborately integrating   Nesterov's extrapolation and  VR techniques, ACC-SADMM provably  achieves a non-ergodic $O(1/K)$ convergence rate which is optimal for  non-smooth problems.  So ACC-SADMM fills the theoretical gap between the stochastic and batch~(deterministic) ADMM.  The original idea to design our ACC-SADMM  is by explicitly considering the snapshot vector $\tx$ into the extrapolation terms.  This is, to some degree,  inspired by~\cite{Katyusha} who proposes an $O(1/K^2)$ stochastic gradient algorithm named Katyusha for convex problems. However, there are  many distinctions between the two algorithms~(please see detailed discussions in Section~\ref{comKrusha}).  Our method is also very efficient in practice since we have  sufficiently considered the noise of  gradient into our acceleration scheme.  For example, we adopt extrapolation as $\y_s^k = \x^k_s+(1-\theta_{1,s}-\theta_2)(\x_s^k - \x_s^{k-1})$ in the inner loop, where $\theta_2$ is a constant and $\theta_{1,s}$  decreases after each whole inner loop, instead of  directly adopting extrapolation as $\y^k = \x^k + \frac{\theta_1^k(1-\theta_1^{k-1})}{\theta_1^{k-1}}(\x^k - \x^{k-1})$ in the original Nesterov's scheme as ~\cite{OPT-SADMM} does. So our extrapolation is more ``conservative'' to tackle the noise of gradients. There are also variants on updating of multiplier and the snapshot vector.  We list the contributions of  our work as follows:
\begin{itemize}
	\item   We propose ACC-SADMM for large scale convex finite-sum problems with linear constraints which  integrates   Nesterov's extrapolation and  VR techniques.  We prove that our algorithm converges in non-ergodic $O(1/K)$ which is optimal for separable linearly constrained non-smooth convex problems.  To our best knowledge, this is the \emph{first} work that achieves a truly accelerated, stochastic convergence rate for constrained convex problems.
	\item  Our algorithm is fast in practice. We have sufficiently considered the noise of gradient into the extrapolation scheme. The memory cost of our method is also low. The experiments  on four bench-mark datasets demonstrate the superiority of our algorithm. We also do experiments on the  Multitask Learning~\cite{multitask} problem  to demonstrate that our algorithm can be used on very large datasets.
\end{itemize}
\vspace{-0.1cm}
\bfseries Notation. \mdseries We denote $\| \x \|$ as the  Euclidean norm, $\langle \x,\y  \rangle = \x^T \y$, $\| \x \|_2=\sqrt{\x^T\x}$,  $\| \x \|_\G= \sqrt{\x^T \G \x}$, and $\langle \x,\y  \rangle_\G = \x^T \G \y$, where $\G\succeq \mathbf{0}$. For a matrix $\X $, $\|\X\|$ is its spectral norm. We use $\I$ to denote the identity matrix. Besides,  a function $f$ has Lipschitz continuous gradients if $\|\nabla f(\x)-\nabla f(\y)\|\leq L\|\x-\y  \|$, which implies \cite{nesterov2013introductory}:
\begin{equation}
f(\y)\leq f(\x )+ \langle \nabla f(\x), \x-\y   \rangle +\frac{L}{2}\|\x-\y\|^2,
\end{equation}
where $\nabla f(\x)$ denotes the gradient of $f$.

\section{Related Works and Preliminary}
\subsection{Accelerated Stochastic Gradient Algorithms}
There are several works in which the authors propose  accelerated, stochastic algorithms for \emph{unconstrained} convex problems. \citet{nitanda2014stochastic} accelerates SVRG~\cite{SVRG} through Nesterov's extrapolation for strongly convex problems. However, their method cannot be extended to general convex problems. Catalyst~\cite{frostig2015regularizing} or APPA~\cite{lin2015universal} reduction also take strategies to obtain faster convergence rate for stochastic convex problems. When the objective function is smooth, these methods can achieve optimal $O(1/K^2)$ convergence rate.  Recently, \citet{Katyusha} and \citet{hien2016accelerated} propose  optimal $O(1/K^2)$ algorithms for general convex problems, named Katyusha and ASMD, respectively. For  $\sigma$-strongly convex problems, Katyusha also meets the optimal $O((n+\sqrt{nL/\sigma})\log \frac{1}{\epsilon}))$ rate.  However,  none of the above algorithms  considers the problems with constraints.

\subsection{Accelerated Batch ADMM Methods}
There are two lines of works which attempt to accelerate Batch ADMM through Nesterov's acceleration schemes. The first line adopts acceleration only on the smooth term~($f_i(\x)$). The works are from ~\cite{ouyang2015accelerated,lu2015fast}.  The convergence rate that they obtain is ergodic $O(\frac{R^2}{K^2}+\frac{D_y}{K})$. The dependence on the smoothness constant is accelerated to $O(1/K^2)$. So these methods are  faster than  ADMM but still remain $O(1/K)$ in the ergodic sense.  The second line is to adopt acceleration on both $f_i(\x)$ and constraints. The resultant algorithm is from \cite{LADM-NE} which is proven  to have a  non-ergodic $O(1/K)$ rate. Since the original ADMM have a tight $O(1/\sqrt{K})$ convergence rate~\cite{ADMM} in the non-ergodic sense, their method is faster theoretically.

\subsection{SADMM and Its Variance Reduction Variants}\label{relaadmm}
We introduce some preliminaries of SADMM. Most SADMM methods alternately minimize the following variant  surrogate of the augmented Lagrangian:
\begin{eqnarray}\label{sLagrangian}
L'(\x_1,\x_2,\olam, \beta) &=&   h_1(\x_1) +\langle \nabla f_1(\x_1), \x_1\rangle  +\frac{L_1}{2}\| \x_1 - \x_1^k\|^2_{\G_1}\\
&& +  h_2(\x_2) +\langle \tna f_2(\x_2), \x_2\rangle  +\frac{L_2}{2}\| \x_2 - \x_2^k\|^2_{\G_2} + \frac{\beta}{2}\| \A_1\x_1+\A_2\x_2-\b + \frac{\olam}{\beta} \|^2, \notag
\end{eqnarray}
where $\tna f_2(\x_2)$ is an estimator of  $\nabla f_2(\x_2)$ from one or a mini-batch of  training samples. So the computation cost for each iteration reduces from $O(n)$ to $O(b)$ instead, where $b$ is the mini-batch size. When $f_i(\x)=0$ and $\G_i=\mathbf{0}$, with $i = 1,2$,  Problem~\eqref{problem} is solved as exact ADMM. When there is no $h_i(\x_i)$,  $\G_i$ is set as the identity matrix $\I$, with $i = 1,2$,  the subproblem in $\x_i$ can be solved through  matrix inversion. This scheme  is  advocated in many SADMM methods~\cite{STOC-ADMM, SAG-ADMM}.  Another common approach is linearization (also called the inexact Uzawa method)~\cite{lin2011linearized,zhang2011unified}, where  $\G_i$ is set as $\eta_i I - \frac{\beta}{L_i} \A^T_i\A_i$ with $\eta_i \geq   1 +\frac{\beta}{L_i} \|\A^T_i\A_i \|$.

For STOC-ADMM~\cite{STOC-ADMM}, $\tna f_2(\x_2)$ is simply set as:
\begin{eqnarray}\label{svrg}
\tna f_2(\x_2) =  \frac{1}{b}\sum_{i_k\in \mathcal{I}_k}\nabla f_{2,i_k}(\x_2),
\end{eqnarray}
where $\mathcal{I}_k$ is the mini-batch of size $b$ from $\{1,2,\cdots,n\}$. 

VR methods~\cite{SDCA-ADMM,SAG-ADMM,SVRG-ADMM} choose more sophisticated gradient estimator to achieve faster convergence rates. As our method bounds the variance through the technique of SVRG~\cite{SVRG},  we introduce SVRG-ADMM~\cite{SVRG-ADMM}, which  uses the gradient estimator as:
\begin{eqnarray}
\tna f_2(\x_2)= \frac{1}{b}\sum_{i_k\in \mathcal{I}_k}\left(\nabla f_{2,i_k}(\x_2) -\nabla f_{2,i_k}(\tx_2) \right)+\nabla f_2(\tx_2),
\end{eqnarray}
where $\tx_2$ is a snapshot vector.  An advantage of SVRG~\cite{SVRG} based methods is its low storage requirement, independent of the number of training samples. This makes them more practical on very large datasets. In our multitask learning experiments, SAG-ADMM~\cite{SAG-ADMM} needs $38.2$TB for storing the weights, and SDCA-ADMM needs $9.6$GB~\cite{SDCA-ADMM} for the dual variables, while the memory cost for our method and SVRG-ADMM is no more than $250$MB.  \citet{SVRG-ADMM} prove that  SVRG-ADMM converges in ergodic $O(1/K)$.  Like batch-ADMM,  in non-ergodic sense, the convergence rate is tight $O(1/\sqrt K)$~(see the discussions in Section \ref{VR-SADMM}).

\section{Our Algorithm}
\subsection{ACC-SADMM}
In this section, we introduce our Algorithm: ACC-SADMM, which is shown in Algorithm~\ref{ACC-SADMM}.  For simplicity, we directly linearize both the smooth term $f_i(\x_i)$ and the augmented term $\frac{\beta}{2}\|\A_1\x_1+\A_2\x_2-\b+\frac{\olam}{\beta}   \|^2$ . It is not hard to extend our method to other schemes mentioned in Section~\ref{relaadmm}.  ACC-SADMM includes two loops. In the inner loop, it updates the primal and dual variables $\x^k_{s,1}$, $\x^k_{s,2}$ and  $\olam^k_s$.  Then in the outer loop, it preserves snapshot vectors $\tx_{s,1}$,  $\tx_{s,2}$ and $\tb_s$, and then resets the initial value of the extrapolation term $\y^0_{s+1}$. Specifically,  in the inner iteration, $\x_1$  is updated as:

\begin{eqnarray}\label{updatax1}
\x^{k+1}_{s,1} &=& \argmin_{x_1} h_1(\x_1)+\langle \nabla f_1(\y^{k}_{s,1}),\x_1\rangle \\ 
&&+ \langle \frac{\beta}{\theta^s_1}\left(\A_1\y_{s,1}^k+\A_2\y_{s,2}^k-\b\right)+\olam^k_s,\A_1^T\x_1  \rangle+\left(\frac{L_1}{2}+\frac{\beta \| \A_1^T\A_1\|}{2\theta^s_1}\right)\|\x_1-\y^{k}_{s,1}\|^2. \notag
\end{eqnarray}
And $\x_2$ is updated using the latest information of $\x_1$, which can be written as:
\begin{eqnarray}\label{updatax2}
\x^{k+1}_{s,2} &=& \argmin_{\x_2} h_2(\x_2)+\langle \tna f_2(\y^{k}_{s,1}),\x_2 \rangle \\ 
&&+ \langle \frac{\beta}{\theta^s_1}\left(\A_1\x_{s,1}^{k+1}+\A_2\y_{s,2}^k-\b\right)+\olam^k_s,\A_2^T\x_2  \rangle +\left(\frac{(1+\frac{1}{b\theta_2})L_2}{2}+\frac{\beta \| \A_2^T\A_2\|}{2\theta^s_1}\right)\|\x_2-\y^{k}_{s,2}\|^2, \notag
\end{eqnarray}
where $\tna f_2(\y_{s,2}^k)$ is obtained by the technique of SVRG~\cite{SVRG} with the form:
\begin{eqnarray}\label{estimate}
\tna f_2(\y_{s,2}^k)=\frac{1}{b}\sum_{i_{k,s}\in \mathcal{I}_{(k,s)}}\left(\nabla f_{2,i_{k,s}}(\y_{s,2}^k- \nabla f_{2,i_{k,s}}(\tx_{s,2}) +\nabla f_2(\tx_{s,2})\right).\notag
\end{eqnarray}
And $\y^{k+1}_s$ is generated as
\begin{eqnarray}\label{geny1}
\y_{s}^{k+1} = \x_{s}^{k+1} +  (1-\theta_{1,s}-\theta_2)(\x_{s}^{k+1} - \x_{s}^{k}),
\end{eqnarray}
when $k\geq0$. One can find that  $1-\theta_{1,s}-\theta_2\leq 1-\theta_{1,s}$.  We do extrapolation in a more ``conservative'' way to tackle  the noise of gradient.  Then the multiplier is updated through Eq.~\eqref{lam1} and Eq.~\eqref{lam2}.  We can find that $\olam^k_s$ additionally  accumulates a compensation term $\frac{\beta\theta_2}{\theta_{1,s}}(\A_1\x_1+\A_2\x_2-\tb_s)$  to ensure $\A_1\x_1+\A_2\x_2$ not to go far from $\A_1\tx_1+\A_2\tx_2$ in the course of iteration.  
\begin{algorithm}[tb]
	\caption{ACC-SADMM}
	\label{ACC-SADMM}
	\begin{algorithmic}
		\STATE  $\!\!\!\!\!\!$ {\textbf{Input:}} epoch length $m>1$, $\beta$,  $\tau=2$, $c=2$, $\x^0_0=\mathbf{0}$, $\tlam^0_0=\mathbf{0}$, $\tx^0=\x^0_0$, $\y^0_0=\x^0_0,$ \\~~~~~~~$\theta_{1,s}=\frac{1}{c+\tau s}$, $\theta_2=\frac{m-\tau}{\tau(m-1)}$. 
		\STATE {\textbf{for}} $s=0$ to $S-1$ {\textbf{do}} 
		\STATE $\quad$ {\textbf{for}} $k=0$ to $m-1$ {\textbf{do}} \vspace{-0.68cm}
		\STATE \begin{eqnarray}\label{lam1}
		\quad\quad\olam^{k}_s = \tlam_s^{k} +\frac{\beta \theta_2}{\theta_{1,s}}\left(\A_1\x^{k}_{s,1}+\A_2\x^{k}_{s,2}-\tb_s\right) .
		\end{eqnarray}\vspace{-0.47cm}
		\STATE $\quad\quad$ Update $\x^{k+1}_{s,1}$ through Eq.~\eqref{updatax1}.
		\STATE$\quad\quad$  Update $\x^{k+1}_{s,2}$ through Eq.~\eqref{updatax2}.  \vspace{-0.6cm}
		\STATE \begin{eqnarray}\label{lam2}
		\quad\quad\tlam_s^{k+1} = \olam_s^{k} + \beta \left(\A_1\x_{s,1}^{k+1}+\A_2\x_{s,2}^{k+1}-\b \right).
		\end{eqnarray} \vspace{-0.6cm}
		\STATE $\quad\quad$ Update $\y^{k+1}_{s}$ through Eq.~\eqref{geny1}.\vspace{0.07cm}
		\STATE $\quad$ {\textbf{end for}} k.
		\STATE  $\quad~\x_{s+1}^0=\x_{s}^m$.\vspace{0.1cm}	
		\STATE $\quad$ Update $\tx_{s+1}$ through Eq.~\eqref{txgen}.	\vspace{-0.6cm}
		\STATE \begin{eqnarray}  \label{lam3}
		\quad~\tlam_{s+1}^{0} = \olam_s^{m-1} + \beta(1-\tau) \left(\A_1\x_{s,1}^{m}+\A_2\x_{s,2}^{m}\!-\b \right).
		\end{eqnarray} \vspace{-0.6cm}
		\STATE $\quad~ \tb_{s+1} =  \A_1 \tx_{s+1,1}  + \A_2 \tx_{s+1,2}.$\vspace{0.07cm}
		\STATE $\quad$ Update $\y^{0}_{s+1}$ through Eq.~\eqref{geny2}.\vspace{0.07cm}
		\STATE {\textbf{end for}} s.
		\STATE $\!\!\!\!\!\!${\textbf{Output:}} 
		\STATE \vspace{-0.9cm} $$\!\!\!\!\!\hat{\x}_S=  \frac{1}{(m-1)(\theta_{1,S}+\theta_2)+ 1}\x^m_S + \frac{\theta_{1,S}+\theta_2}{(\!m-\!1)(\theta_{1,S}+\theta_2)+1} \sum_{k=1}^{m-1}\x^k_S.$$\vspace{-0.3cm}
	\end{algorithmic}
\end{algorithm}

In the outer loop, we set the snapshot vector $\tx_{s+1}$  as:
\begin{eqnarray}\label{txgen}
	\tx_{s+1}=  \frac{1}{m}\left(\left[1-\frac{(\tau-1)\theta_{1,{s+1}}}{\theta_2} \right]\x^m_s+\left[1+\frac{(\tau-1)\theta_{1,{s+1}}}{(m-1)\theta_2}\right]\sum_{k=1}^{m-1}\x^k_s\right).
\end{eqnarray}
$\tx_{s+1}$ is not the average of $\{\x^k_s \}$,  different from most  SVRG-based methods~\cite{SVRG,SVRG-ADMM}. The way of generating  $\tx$ guarantees a faster convergence rate for the constraints.  Then  at the last step, we reset  $\y_{s+1}^0$ as:
\begin{eqnarray}\label{geny2}
\y_{s+1}^0=(1-\theta_2)\x^m_{s} +\theta_2 \tx_{s+1}+\frac{\theta_{1,s+1}}{\theta_{1,s}} 
\left[(1-\theta_{1,s})\x^m_{s}-(1-\theta_{1,s}-\theta_2)\x^{m-1}_{s} -\theta_2 \tx_s   \right].\notag
\end{eqnarray}
The whole algorithm is shown Algorithm~\ref{ACC-SADMM}.

\subsection{Intuition}
Though Algorithm \ref{ACC-SADMM} is a little complex at the first sight, our intuition to design the algorithm is straightforward.  To bound the variance, we use the technique of SVRG~\cite{SVRG}. Like \cite{SVRG,Katyusha}, the variance of gradient is bounded through 
\begin{eqnarray}\label{variance}
 \E_{i_{k}}\left(\|\nabla f_2(\y_2^k) - \tna f_2(\y_2^k)\|^2\right)\leq \frac{2L_2}{b} \left[ f_2(\tx_2) -  f_2(\y^k_2) +   \langle \nabla  f_2(\y_2^k), \y^k_2-\tx_2  \rangle \right], 
\end{eqnarray}
where $\E_{i_k}$ indicates that the expectation is  taken over the random choice of $i_{k,s}$,  under the condition that $\y_2^k$, $\tx_2$ and $\x^k_2$ (the  randomness in the first $sm+k$ iterations are fixed) are known.
We first consider the case that there is no linear constraint. Then by the convexity of $F_1(\x_1)$~\cite{FISTA}, we have
\begin{eqnarray}
F_1(\x^{k+1}_1)\leq F_1(\uu_1) + \frac{L_1}{2}\| \x^{k+1}_1-\y^{k+1}_1   \|^2-L_1\langle\x^{k+1}_1-\y^{k}_1, \x^{k+1}_1-\uu_1\rangle, 
\end{eqnarray}
Setting $\uu_1$  be $\x_1^{k}$, $\tx_1$ and $\x_1^*$, respectively, then multiplying the three inequalities by $(1-\theta_1-\theta_2)$, $\theta_2$, and $\theta_1$, respectively,  and adding them,  we have
\begin{eqnarray}\label{FF1}
F_1(\x^{k+1}_1)&\leq& (1-\theta_1-\theta_2)F_1(\x^k_1)+ \theta_2 F_1(\tx_1)+\theta_1 F_1(\x_1^*) \notag\\
&& \!\!\!\!\!\!\!\!\! -L_1 \langle  \x^{k+1}_1\!-\!\y^{k}_1, \x^{k+1}_1\!-\!(1-\theta_1-\theta_2)\x^k_1-\theta_2 \tx_1\!-\!\theta_1\x_1^*\rangle+\frac{L_1}{2}\|\x^{k+1}_1-\y^k_1  \|^2.
\end{eqnarray}
where $\theta_1$ and $\theta_2$ are  undetermined coefficients. Comparing with Eq.~\eqref{variance},  we can find  that there is one more term $\langle \nabla f_2(\y_2^k), \y_2^k - \tx_2\rangle$ that we need to eliminate.  To solve this issue,  we  analyse the points at $\w^k = \y_2^k +\theta_3 (\y_2^k-\tx_2)$ and $\z^{k+1}=\x_2^{k+1}+\theta_3 (\y_2^k-\tx_2)$, where $\theta_3$ is an undetermined coefficient. When $\theta_3>0$,  $\w^k$ and $\z^{k+1}$ is  closer to $\tx_2$ compared with $\y_2^k$ and $\x_2^{k+1}$.  Then by the convexity of $F_2(\x_2)$, we can generate a negative $\langle \nabla f_2(\y_2^k), \y_2^k - \tx_2\rangle$, which can help to eliminate the variance term. 
Next we consider the multiplier term. To construct a recursive term of $L(\x^{k+1}_{s,1},\x^{k+1}_{s,2},\olam^*) - (1-\theta_{1,{s}}-\theta_2) L(\x^{k}_{s,1},\x^{k}_{s,2},\olam^*) -\theta_2L(\tx_{s,1},\tx_{s,2},\olam^*)$, where $L(\x_1,\x_2,\olam)$ satisfies Eq.~\eqref{lam6}, the multiplier should satisfy the following equations:
\begin{eqnarray} \label{lam4}
\hlam^{k+1}_s-\hlam^{k}_s= \frac{\beta }{\theta_{1,s}}\A\left(\x^{k+1}_s\!-\!(1-\theta_1-\theta_2)\x^{k}_s-\theta_1\x^* -\theta_2\tx_s\right),\notag
\end{eqnarray}
and
\begin{eqnarray}\label{lam5}
\hlam^{k+1}_s = \olam^k_s +\frac{\beta}{\theta_{1,s}}(\A\x^{k+1}_s-\b),
\end{eqnarray}
where $\hlam^{k}_s$ is undetermined and
\begin{eqnarray}\label{lam6}
L(\x_1,\x_2,\olam) =  F(\x_1,\x_2) + \langle \olam,\A_1\x_1+\A_2\x_2-\b \rangle,
\end{eqnarray}
is the Lagrangian function.  By introducing a new variable $\tlam^k_s$, then setting 
\begin{eqnarray}\label{lam7}
\hlam^{k}_s = \tlam^{k}_s +\frac{\beta(1-\theta_{1,s})}{\theta_{1,s}}(\A\x^k_s-\b),
\end{eqnarray}
and with Eq.~\eqref{lam1} and  Eq.~\eqref{lam2},  Eq.~\eqref{lam4} and  Eq.~\eqref{lam5} are satisfied. Then Eq.~\eqref{lam3} is obtained as we need $\hlam^{0}_s=\hlam^{m}_{s-1}$ when $s\geq 1$.

\section{Convergence Analysis}
In this section, we give the convergence results of ACC-SADMM. The proof can be found in Supplementary Material.  We first analyse each inner iteration. The result is shown in Lemma~\ref{tolist}, which connects $\x^{k}_s$ to $\x^{k+1}_s$.
\begin{lemma}\label{tolist}
	Assume that $f_{1}(\x_1)$ and  $f_{2,i}(\x_2)$ with $i\in\{1,2,\cdots,n\}$  are convex and have Lipschitz continuous gradients.  $L_1$ is the Lipschitz constant of $f_{1}(\x_1)$.   $L_2$ is the Lipschitz constant of $f_{2,i}(\x_2)$ with $i\in\{1,2,\cdots,n\}$ . $h_1(\x_1)$ and $h_2(\x_2)$ is also  convex. For Algorithm~\ref{ACC-SADMM}, in any epoch, we have
	\begin{eqnarray}
	&&\E_{i_k}\left[L(\x^{k+1}_1,\x^{k+1}_2,\olam^*)\right] - \theta_2 L(\tx_1,\tx_2,\olam^*) -(1-\theta_2 - \theta_1)L(\x^{k}_1,\x^{k}_2,\olam^*)\notag\\ 
	&\leq& \frac{\theta_1}{2\beta}\left(\|\hlam^k-\olam^*\|^2- \E_{i_k}\left[\|\hlam^{k+1}-\olam^*\|^2\right] \right)\notag\\
	&& +\frac{1}{2}\|\y_1^{k}-(1-\theta_1-\theta_2)\x^k_1-\theta_2\tx_1-\theta_1\x_1^*\|^2_{\G_3}-\frac{1}{2}\E_{i_k}\left(\|\x_1^{k+1}-(1-\theta_1-\theta_2)\x^k_1-\theta_2\tx_1-\theta_1\x_1^*\|^2_{\G_3}\right)\notag\\
	&&+ \frac{1}{2} \|\y_2^{k}-(1-\theta_1-\theta_2)\x^k_2-\theta_2\tx_2-\theta_1\x_2^*\|^2_{\G_4}-\frac{1}{2} \E_{i_k} \left(\|\x_2^{k+1}-(1-\theta_1-\theta_2)\x^k_2-\theta_2\tx_2-\theta_1\x_2^*\|^2_{\G_4}\right),\notag
	\end{eqnarray}
	where $L(\x_1,\x_2,\olam)$ satisfies Eq.\eqref{lam6} and $\hlam$  satisfies Eq.\eqref{lam7}, $\G_3=\left(L_1+\frac{\beta\| \A_1^T \A_1\|}{\theta_1}\right)\I-\frac{\beta\A_1^T\A_1}{\theta_1}$,  and  $\G_4=\left((1+\frac{1}{b\theta_2}) L_2+\frac{\beta\| \A_2^T \A_2\|}{\theta_1}\right)\I$. We have  ignored subscript $s$ as $s$ is fixed in each epoch.
\end{lemma}
Then Theorem~\ref{convergence1} analyses ACC-SADMM in the whole iteration, which is the key convergence result of the paper.
\begin{theorem}\label{convergence1}
	If the conditions in Lemma~\ref{tolist} hold, then we  have
	\begin{eqnarray}\label{final1}
	&&\E\left(\frac{1}{2\beta}\|\frac{\beta m}{\theta_{1,{S}}}\left(\A\hat{\x}_{S}\!-\!\b\right)-\frac{\beta(m\!-\!1)\theta_{2}}{\theta_{1,{0}}}\left(\A\x^0_{0}-\b\right) +\tlam^0_0-\olam^*\! \|^2 +\frac{m}{\theta_{1,{S}}}\left(F(\hat{\x}_{S})-F(\x^*)  +\langle \olam^*, \A\hat{\x}_S -\b\rangle\right)\right) \notag\\
	&\leq& C_3\left(F(\x_{0}^0)-F(\x^*)  +\langle \olam^*, \A\x^0_0 -\b\rangle\right)+\frac{1}{2\beta}\|\tlam^0_0 +\frac{\beta(1-\theta_{1,{0}})}{\theta_{1,{0}}}(\A\x_0^0-\b) -\olam^*  \|^2\notag\\
	&&+\frac{1}{2}\|\x^0_{0,1}-\x^*_1   \|^2_{\left(\theta_{1,0}L_1+\| \A_1^T \A_1\|\right)\I-\A_1^T\A_1}+\frac{1}{2}\| \x^0_{0,2}-\x^*_2  \|^2_{\left((1+\frac{1}{b\theta_2}) \theta_{1,0}L_2+\| \A_2^T \A_2\|\right)\I},
	\end{eqnarray}
	where $C_3= \frac{1-\theta_{1,0}+(m-1)\theta_2}{\theta_{1,0}}$.
\end{theorem}
Corollary~\ref{convergence2} directly demonstrates that ACC-SADMM have a non-ergodic $O(1/K)$ convergence rate.
\begin{corollary}\label{convergence2}
	If the conditions in Lemma~\ref{tolist} holds,  we 
	have 
	\begin{eqnarray}
	\E\|F(\hat{\x}_{S})-F(\x^*)| &\leq& O(\frac{1}{S}),\notag\\ 
	\E\| \A\hat{\x}_{S}-\b\| &\leq& O(\frac{1}{S}).
	\end{eqnarray}
\end{corollary}
We can find that $\hat{\x}_S$ depends on the latest $m$ information of $\x^k_S$.  So our convergence results is in non-ergodic sense, while the analysis for SVRG-ADMM~\cite{SVRG-ADMM} and SAG-ADMM~\cite{SAG-ADMM} is in ergodic sense, since they consider the point $\hat{\x}_{S} = \frac{1}{mS}\sum_{s=1}^S\sum_{k=1}^m \x^k_s$, which is the convex combination of $\x^k_s$ over \emph{all} the  iterations. 

Now we directly use the theoretical results of~\cite{LADM-NE} to demonstrate that our algorithm is optimal when there exists non-smooth term in the objective function.
\begin{theorem}\label{convergence3}
	For the following problem:
	\begin{eqnarray}\label{problem2}
	\min_{\x_1,\x_2}F_1(\x_1) +F_2(\x_2), ~~ \text{s.t.}~~\x_1-\x_2=\mathbf{0},
	\end{eqnarray}
	let the ADMM type algorithm to solve it be:
	\begin{itemize}
		\item  Generate $\olam^k_2$ and $\y^k_2$  in any way, 
		\item $\x^{k+1}_1=\text{Prox}_{F_1/\beta^k}\left(\y^k_2-\frac{\olam^k_2}{\beta^k} \right),$
		\item  Generate $\olam^{k+1}_1$ and $\y^{k+1}_1$ in any way,
		\item $\x^{k+1}_2=\text{Prox}_{F_2/\beta^k}\left(\y^{k+1}_1-\frac{\olam^{k+1}_1}{\beta^k} \right).$
	\end{itemize}
	Then there exist convex functions $F_1$ and $F_2$ defined on $\mathcal{X}=\{\x\in \mathcal{R}^{6k+5} : \|\x\|\leq B\}$ for the above general ADMM method, satsifying
	\begin{eqnarray}
L\| \hat{\x}^k_2- \hat{\x}^k_1   \| + | F_1(\hat{\x}^k_1)-F_1(\x^*_1) +F_1(\hat{\x}^k_2)-F_2(\x^*_2) |\geq \frac{LB}{8(k+1)},
	\end{eqnarray}
	where $ \hat{\x}^k_1 = \sum_{i=1}^k \alpha^i_1\x^i_1$ and  $ \hat{\x}^k_2 = \sum_{i=1}^k \alpha^i_2\x^i_2$ for any $\alpha^i_1$ and $\alpha^i_2$ with $i$  from $1$ to $k$.
\end{theorem}
Theorem~\ref{convergence3} is  Theorem~$11$ in~\cite{LADM-NE}. More details can be found in it. Problem~\eqref{problem2} is a special case of Problem~\eqref{problem} as we can set each $F_{2,i}(\x_2)=F(\x_2)$ with $i = 1,\cdots,n$ or set $n=1$.  So there is no better ADMM type algorithm which converges faster than $O(1/K)$ for Problem~\eqref{problem}.
\section{Discussions}
We discuss some properties of ACC-SADMM and make further comparisons with some related methods.
\subsection{Comparison with Katyusha}\label{comKrusha}
As we have mentioned in Introduction,  some intuition of our algorithm is inspired by  Katyusha~\cite{Katyusha}, which obtains an $O(1/K^2)$ algorithm for convex finite-sum problems.  However, Katyusha cannot solve the problem with linear constraints.  Besides, Katyusha uses the  Nesterov's second scheme to accelerate the algorithm while our method conducts acceleration  through  Nesterov's extrapolation~(Nesterov's first scheme).   And our proof uses the technique of \cite{Tseng}, which is different from \cite{Katyusha}.    Our algorithm can be easily extended to unconstrained convex finite-sum  and can also obtain a $O(1/K^2)$ rate but belongs to the  Nesterov's first scheme~\footnote{We follow ~\cite{Tseng} to name the  extrapolation scheme as Nesterov's first scheme and the three-step scheme~\cite{nesterov1988approach} as the Nesterov's second scheme.}.
\subsection{The Importance of Non-ergodic $O(1/K)$}\label{VR-SADMM}
\begin{figure}[t]
	\centering 
	\subfigure[a9a-training]{ 
		{\label{fig:1}} 
		\includegraphics[width=3.4cm,height=3cm]{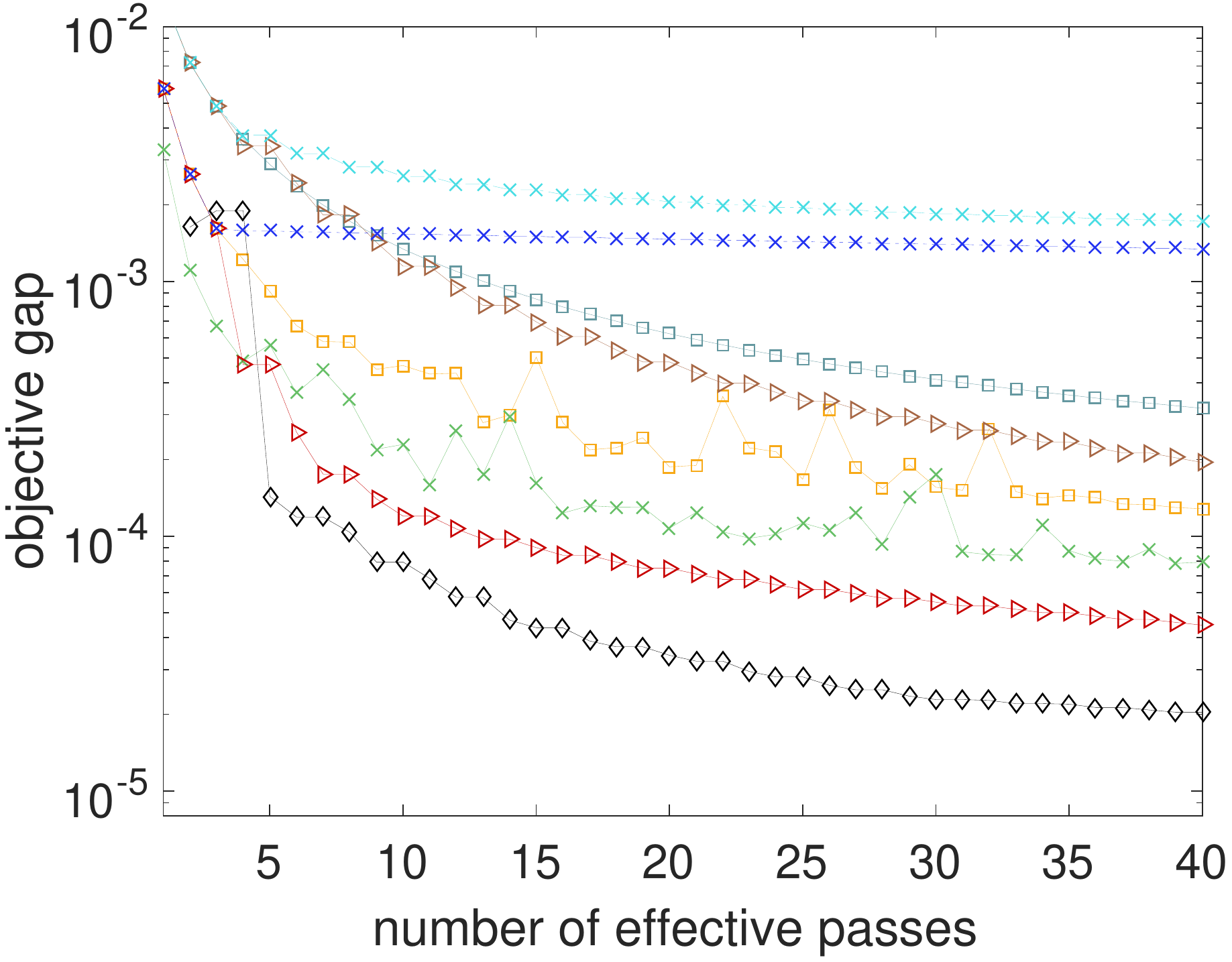}} 
	\subfigure[covertype-training]{ 
		\label{fig:2} 
		\includegraphics[width=3.4cm,height=3cm]{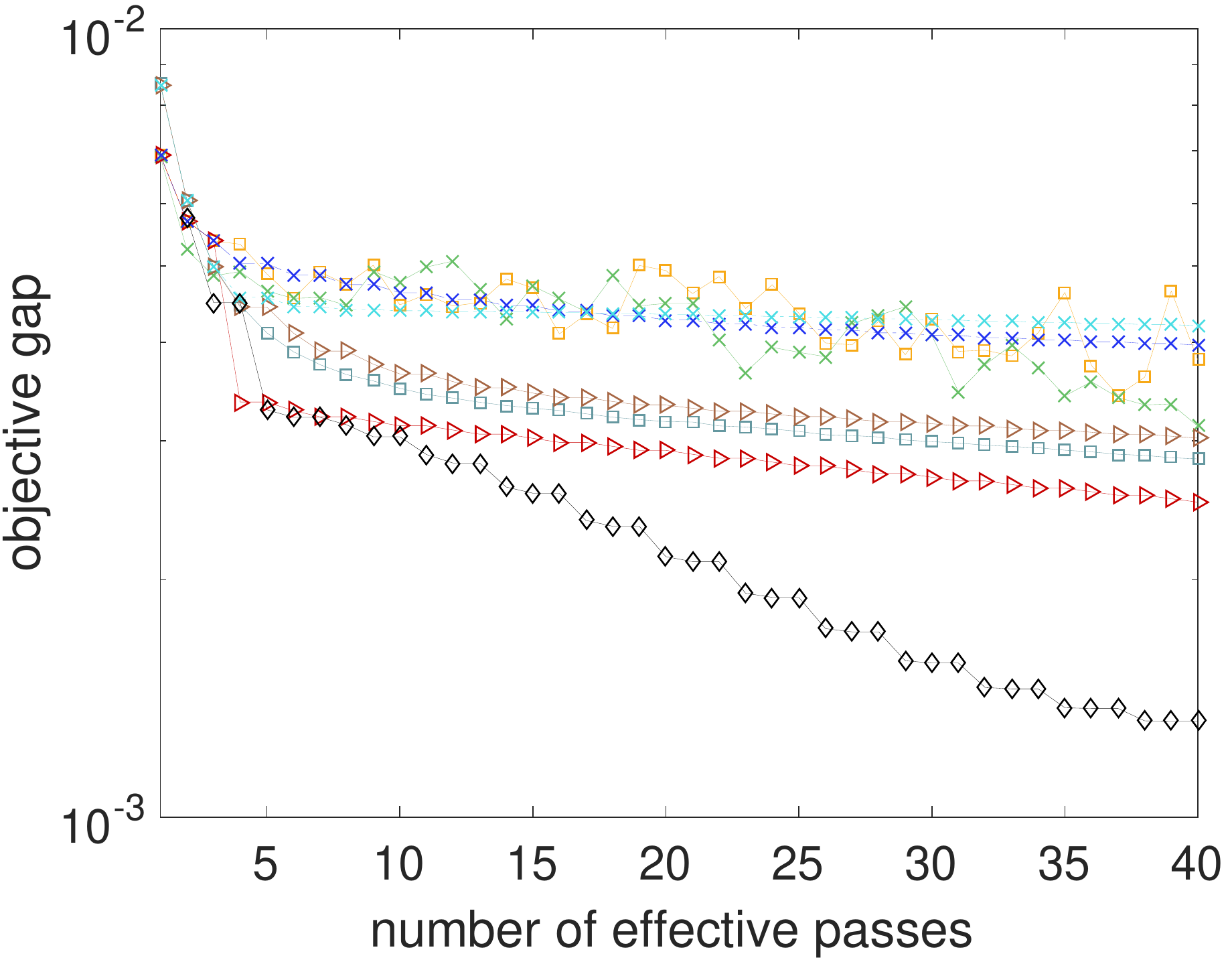}  
	}    
	\subfigure[mnist-training]{ 
		\label{fig:3} 
		\includegraphics[width=3.4cm,height=3cm]{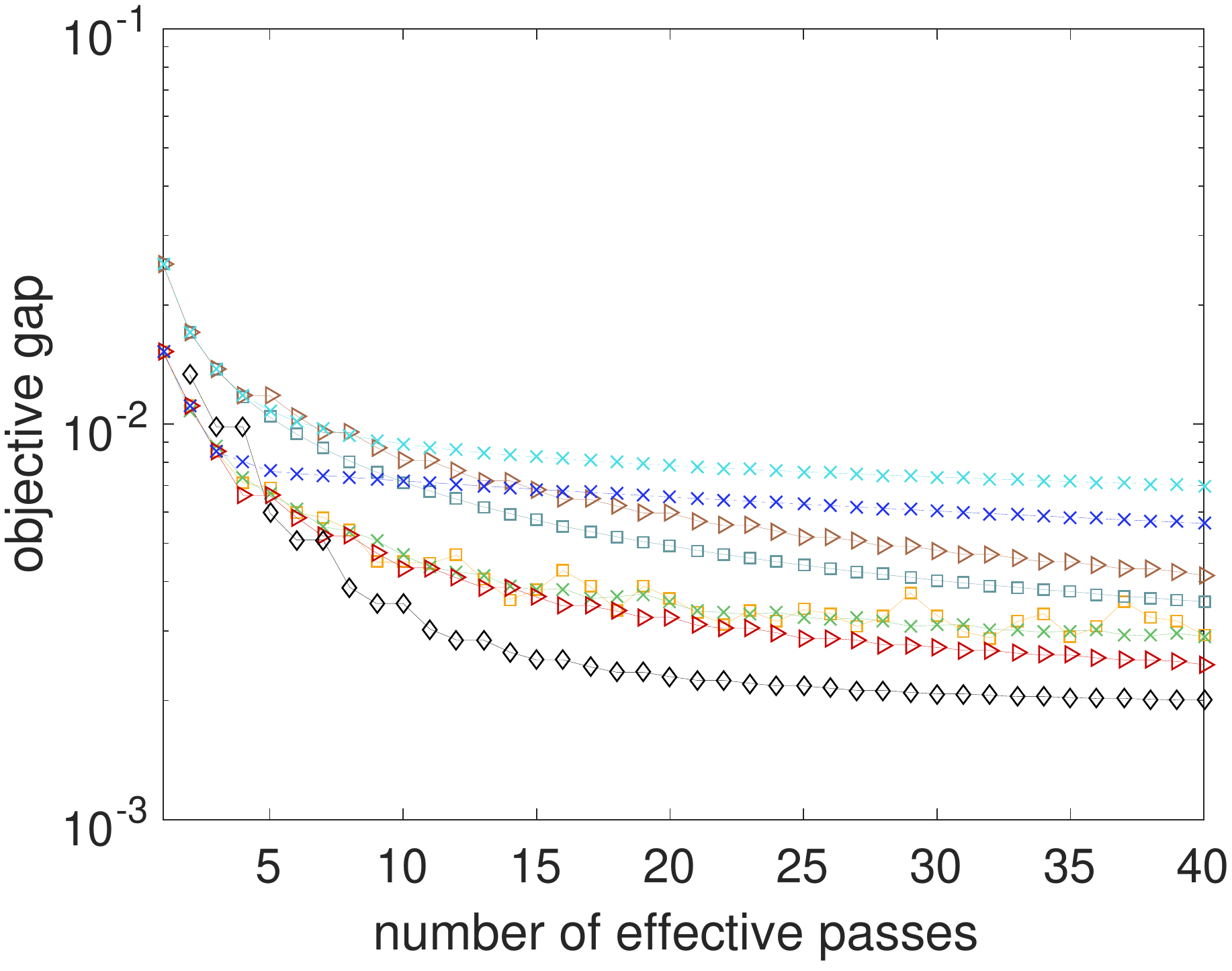}}  
	\subfigure[dna-training]{ 
		\label{fig:4} 
		\includegraphics[width=3.4cm,height=3cm]{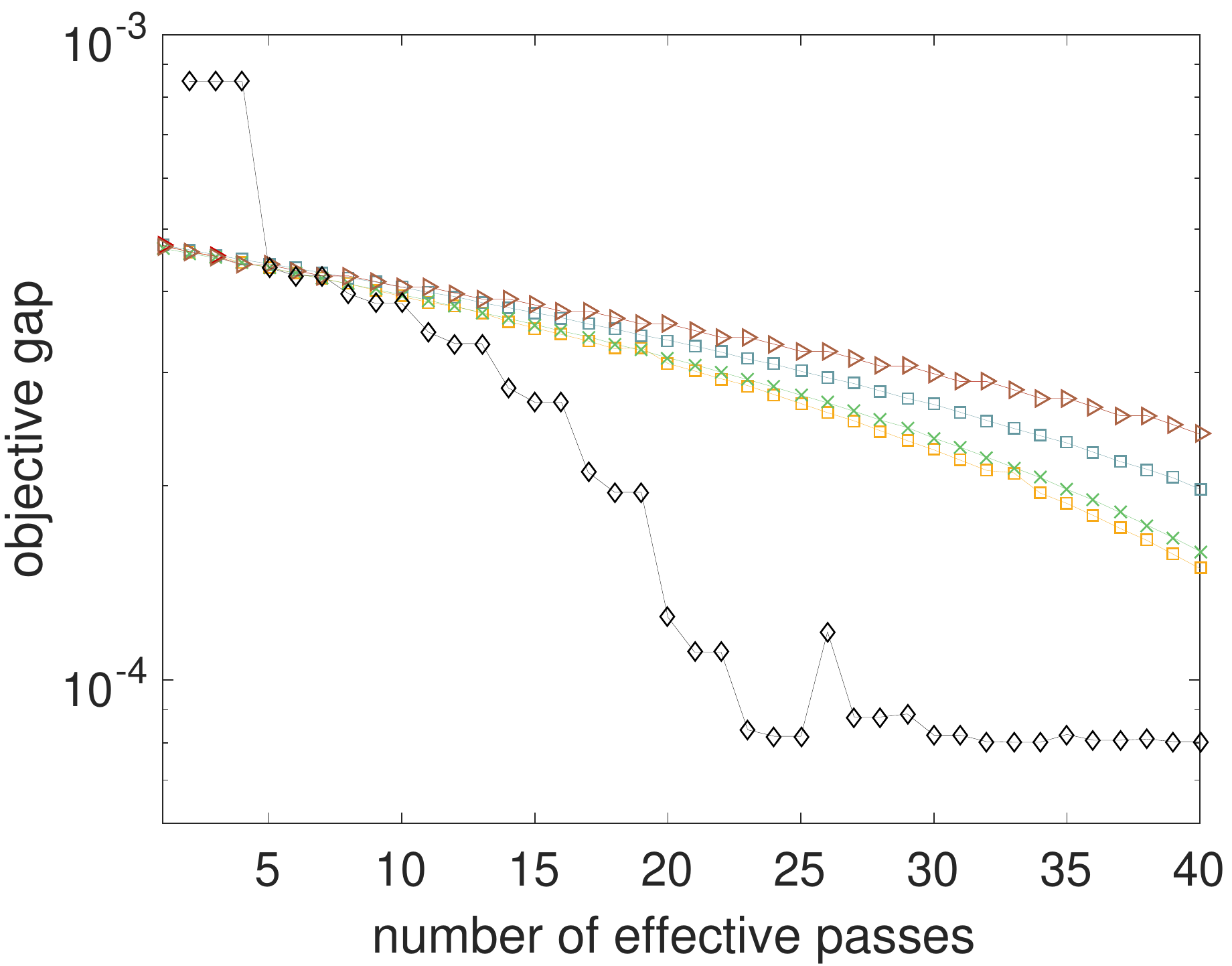}}   
\vspace{-0.3cm}
	\subfigure[a9a-testing]{ 
		{\label{fig:5}} 
		\includegraphics[width=3.4cm,height=3cm]{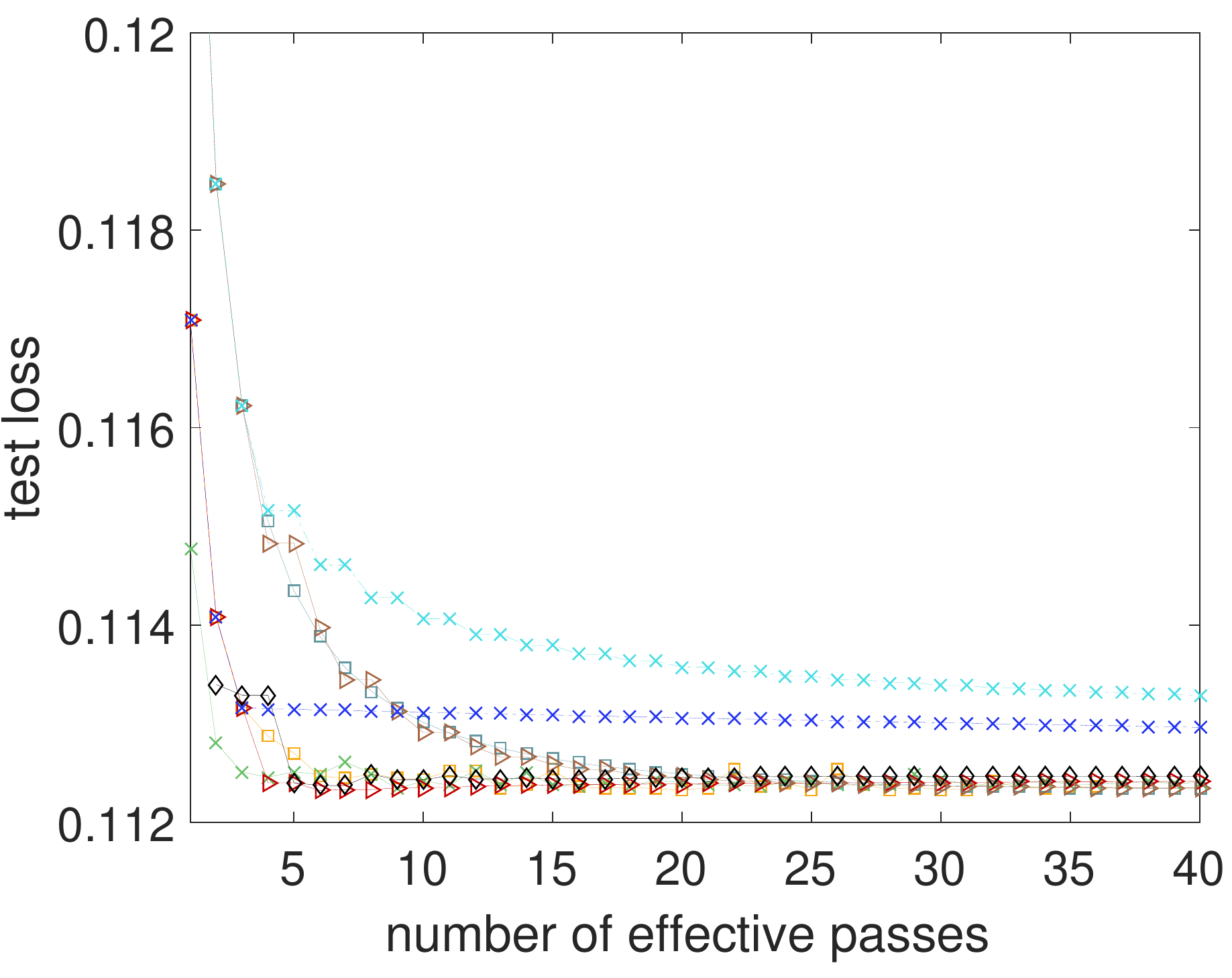}} 
	\subfigure[covertype-testing]{ 
		\label{fig:6} 
		\includegraphics[width=3.4cm,height=3cm]{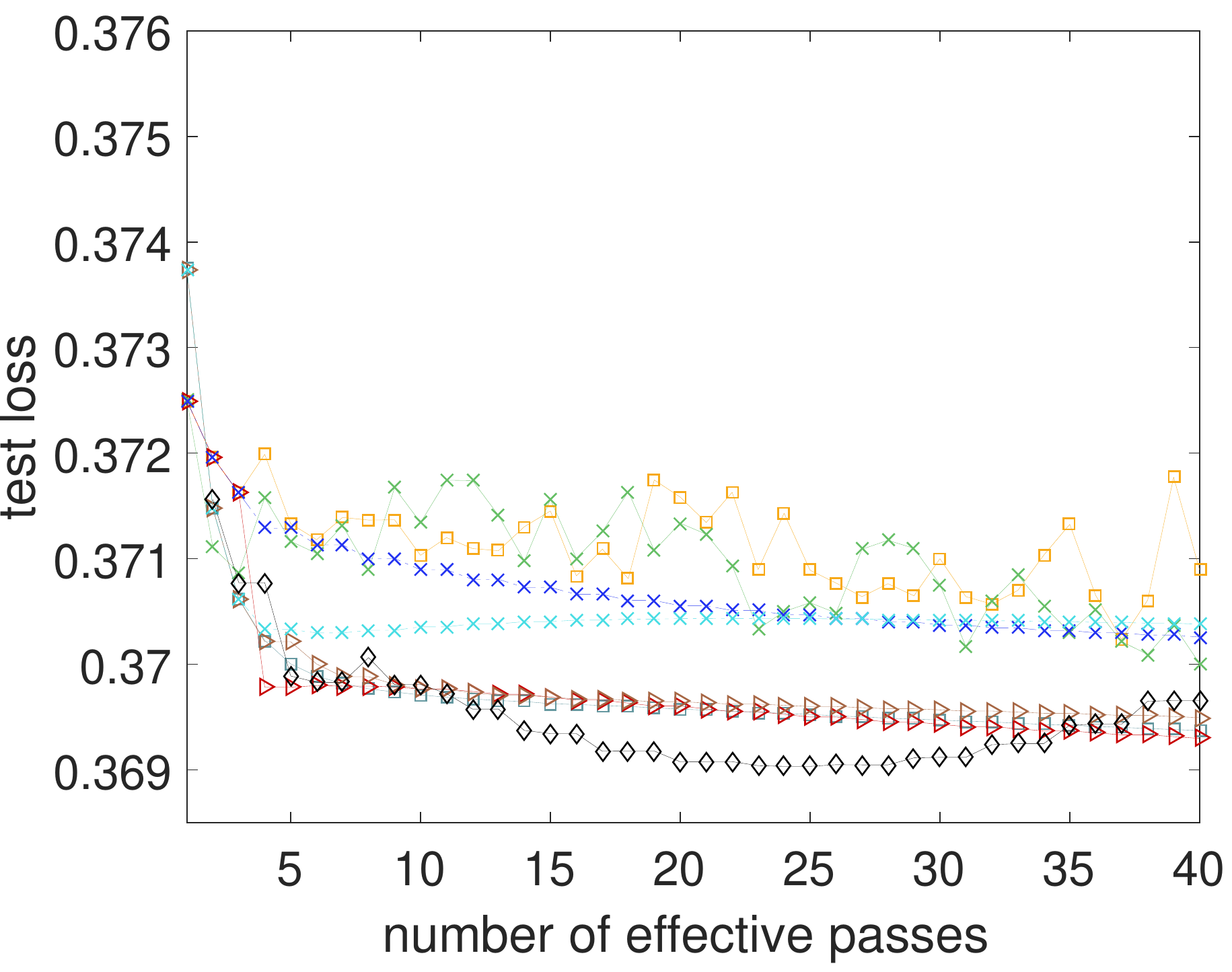}  
	}    
	\subfigure[mnist-testing]{ 
		\label{fig:7} 
		\includegraphics[width=3.4cm,height=3cm]{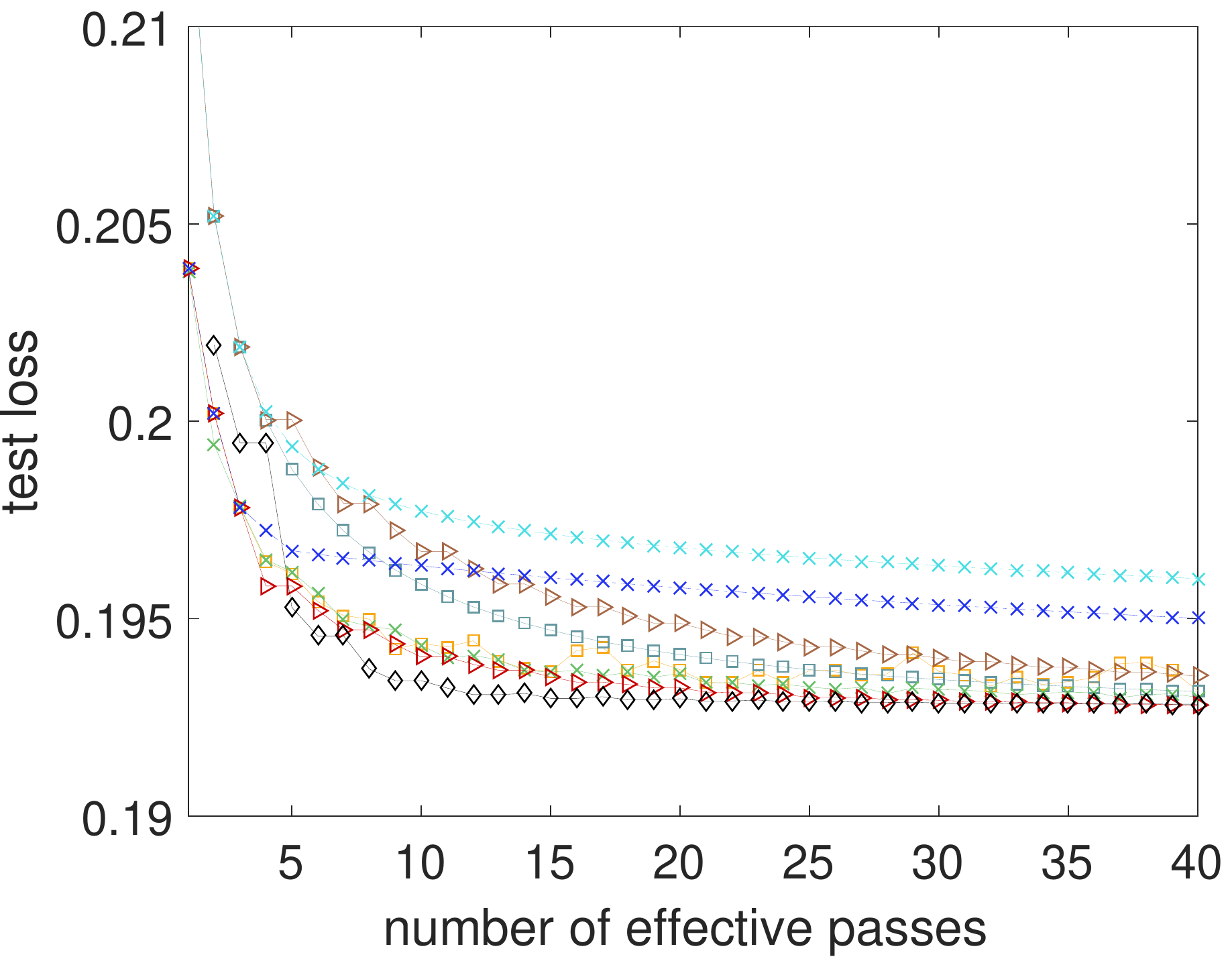}}  
	\subfigure[dna-testing]{ 
		\label{fig:8} 
		\includegraphics[width=3.4cm,height=3cm]{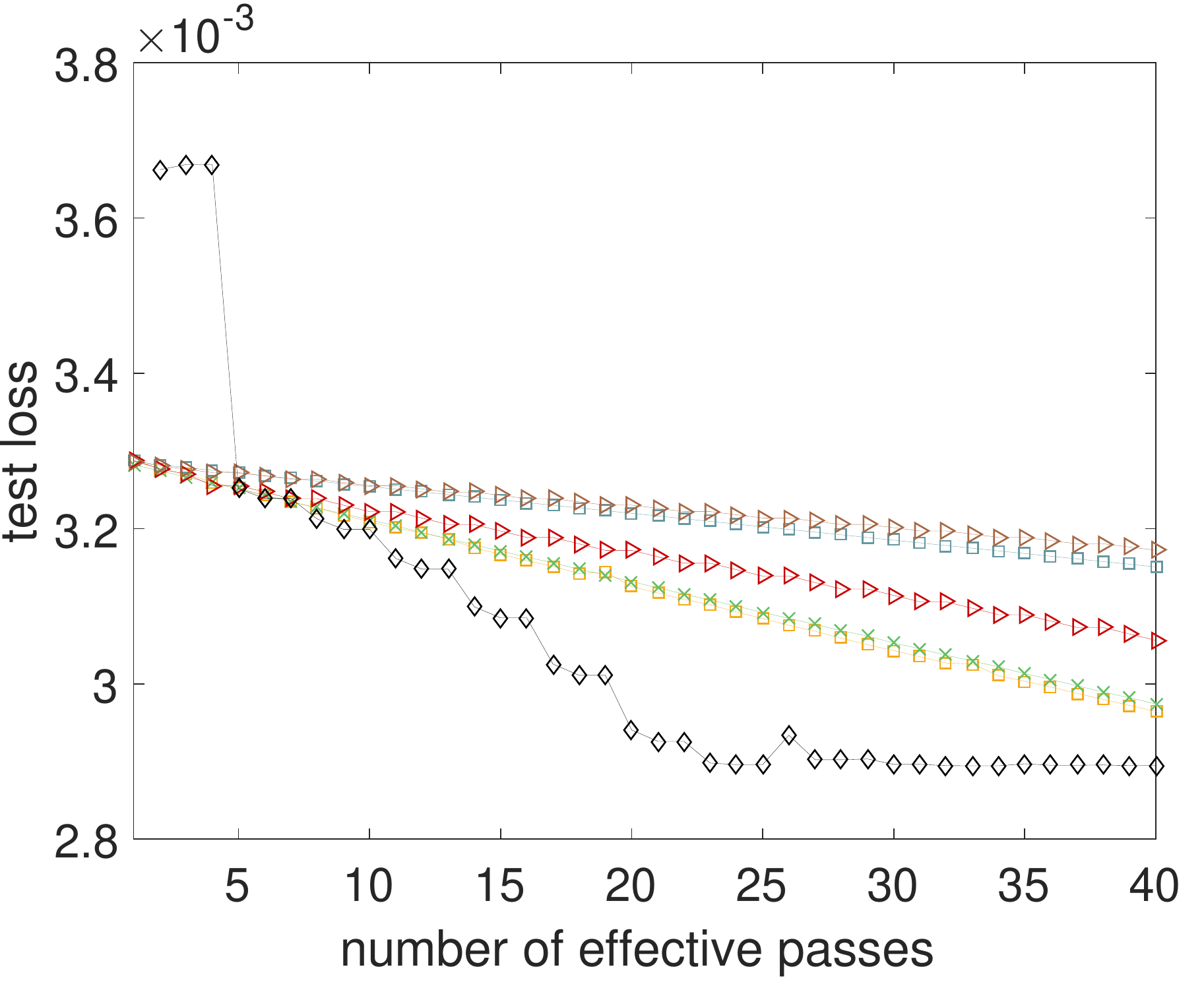}}
	\subfigure{
		\label{fig:9} 
		\includegraphics[width=0.9\linewidth]{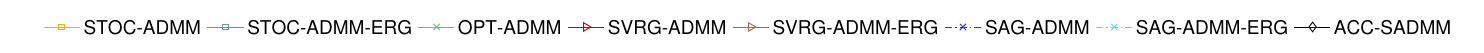}}  
	\vspace{-0.1in}   
	\caption{Experimental results of solving  the original Lasso problem (Eq.~\eqref{lasso})~(Top: objective gap; Bottom: testing loss). The step sizes are set based on the theoretical guidances. The computation time has included the cost of calculating full gradients for SVRG based methods. SVRG-ADMM and SAG-ADMM are initialized by running STOC-ADMM for $\frac{3n}{b}$ iterations. ``-ERG'' represents the ergodic results for the corresponding algorithms.}
	\label{original lasso} 
\end{figure}

\begin{figure}[ht]
	\centering 
	\subfigure[a9a-training]{ 
		{\label{fig:1}} 
		\includegraphics[width=3.4cm,height=3cm]{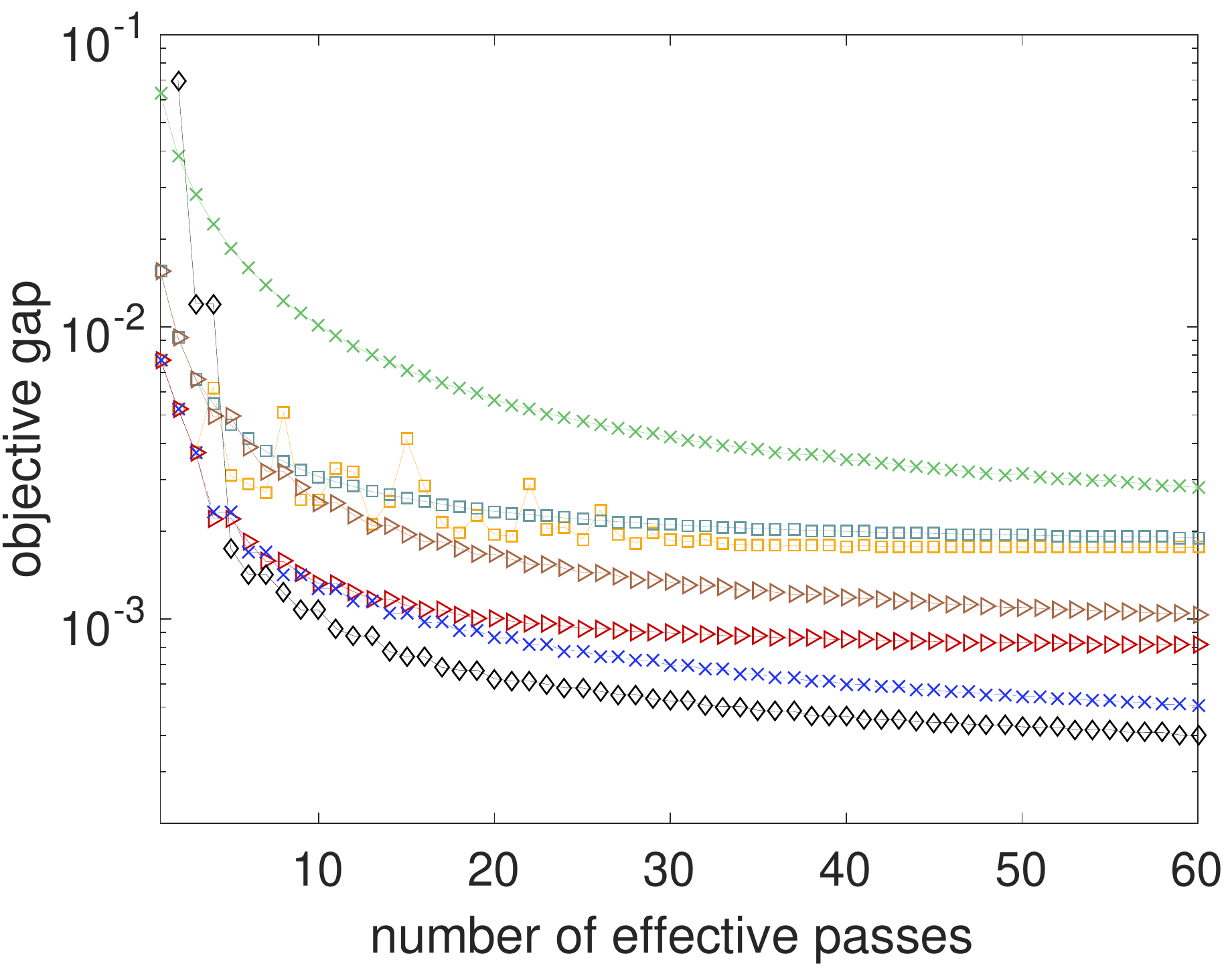}} 
	\subfigure[covertype-training]{ 
		\label{fig:2} 
		\includegraphics[width=3.4cm,height=3cm]{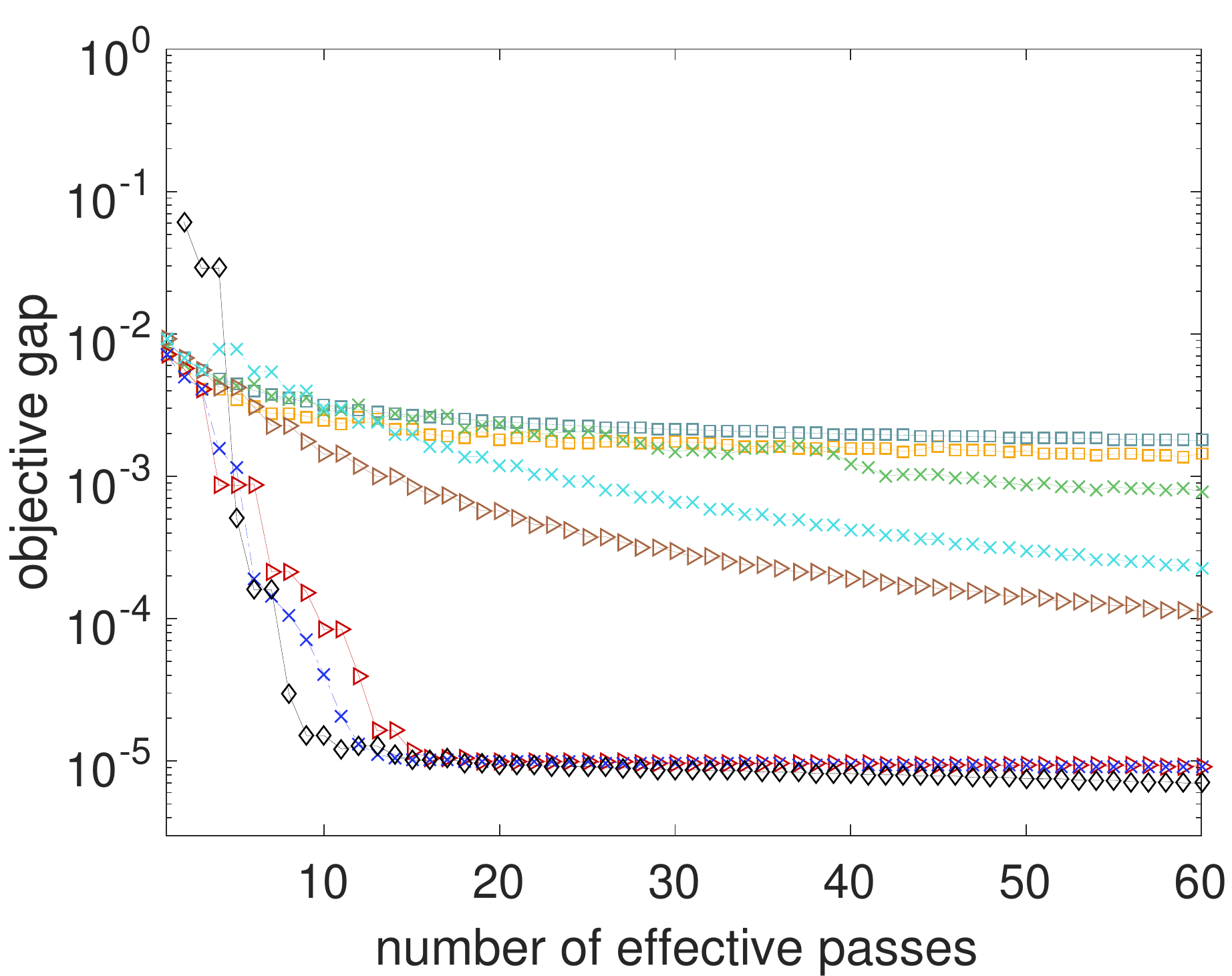}  
	}    
	\subfigure[mnist-training]{ 
		\label{fig:3} 
		\includegraphics[width=3.4cm,height=3cm]{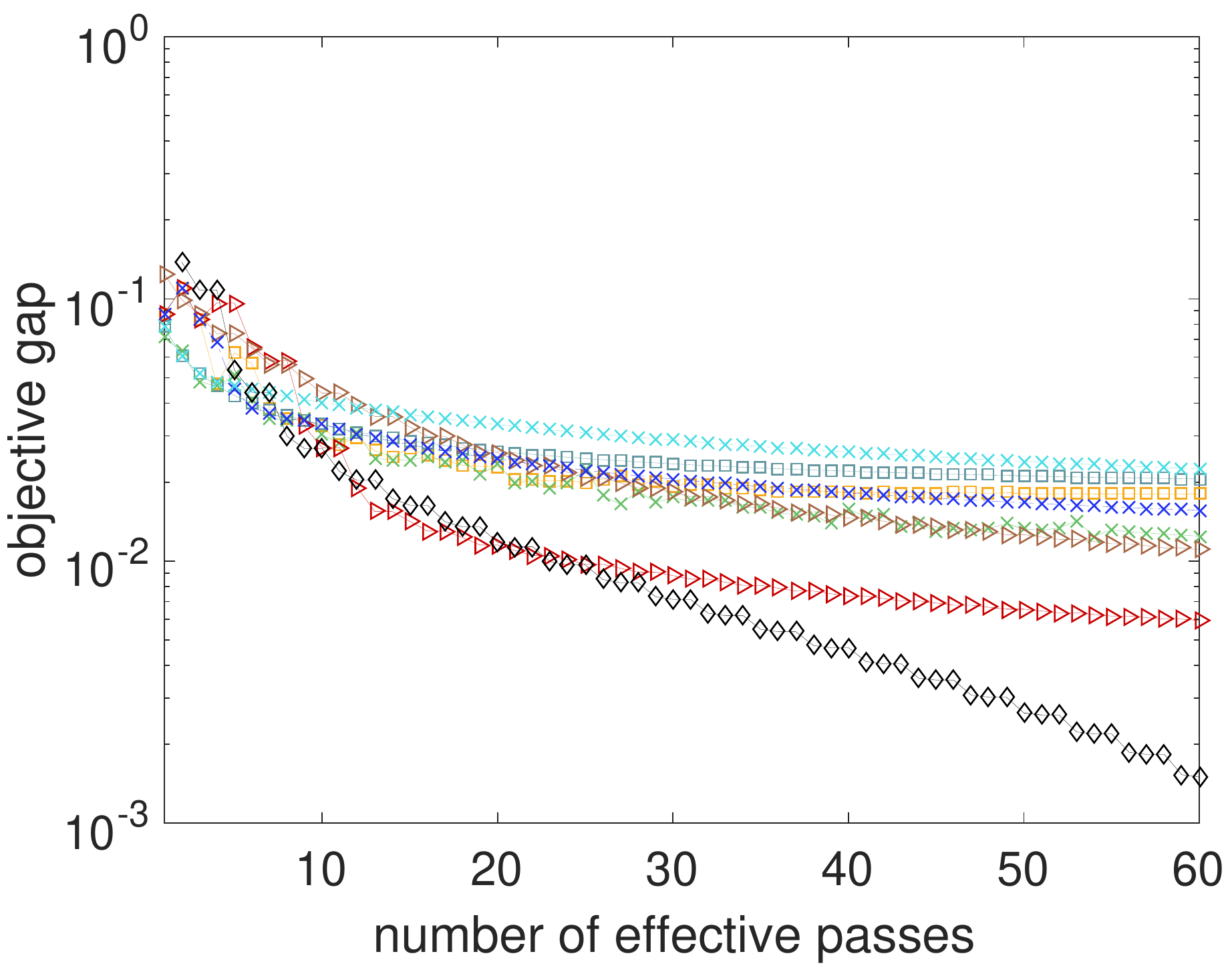}}  
	\subfigure[dna-training]{ 
		\label{fig:4} 
		\includegraphics[width=3.4cm,height=3cm]{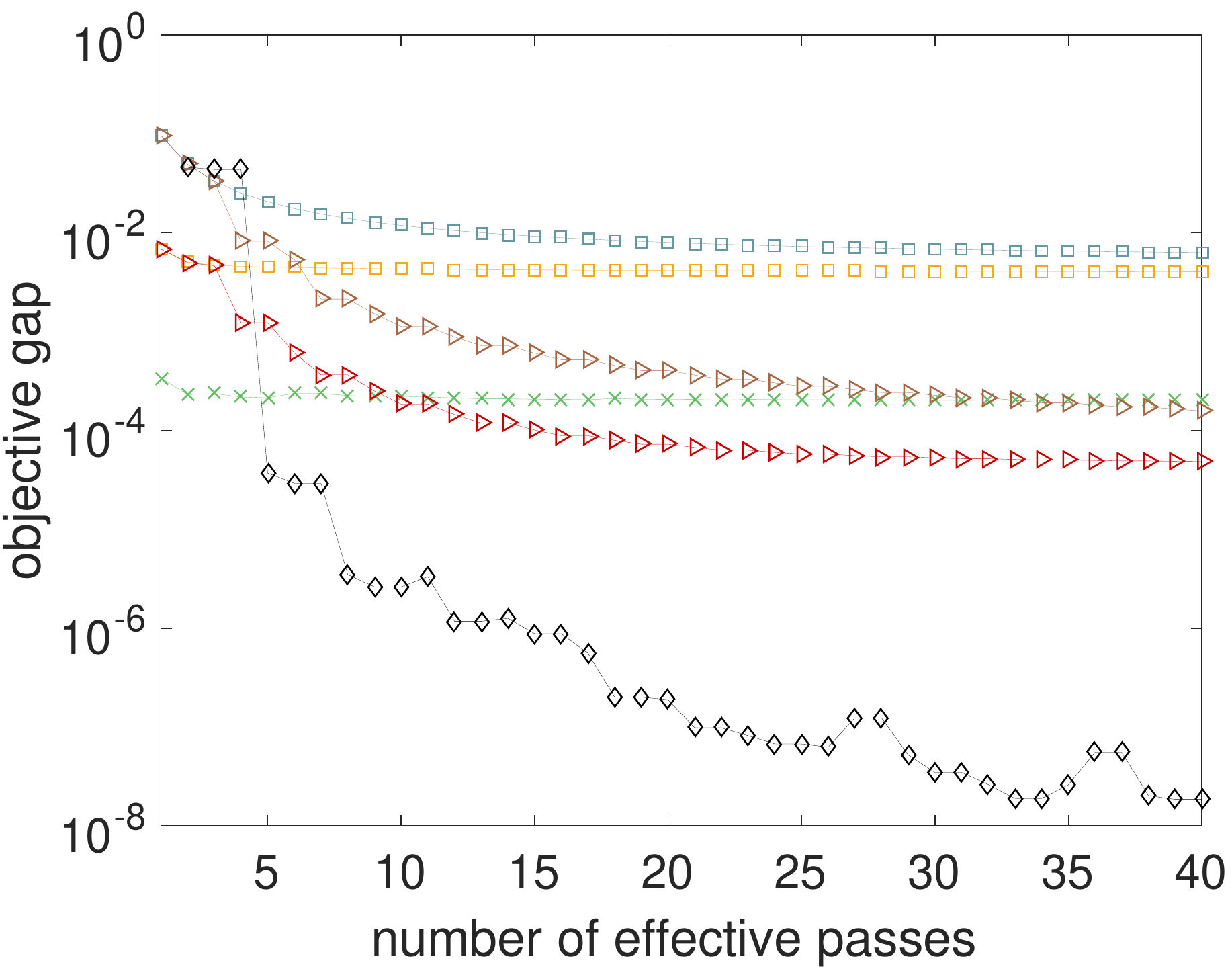}}     
	\vspace{-0.3cm}

	\subfigure[a9a-testing]{ 
		{\label{fig:5}} 
		\includegraphics[width=3.4cm,height=3cm]{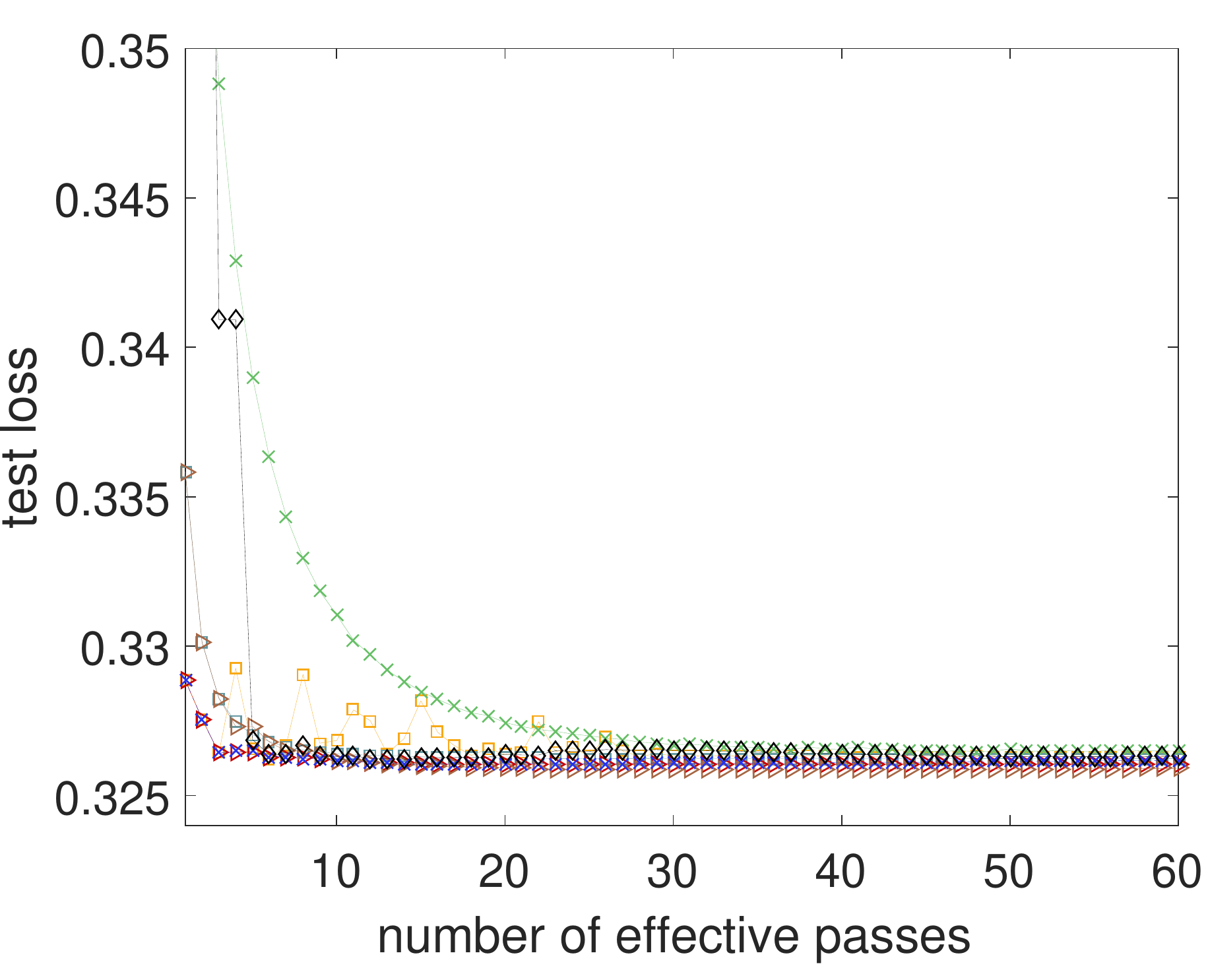}} 
	\subfigure[covertype-testing]{ 
		\label{fig:6} 
		\includegraphics[width=3.4cm,height=3cm]{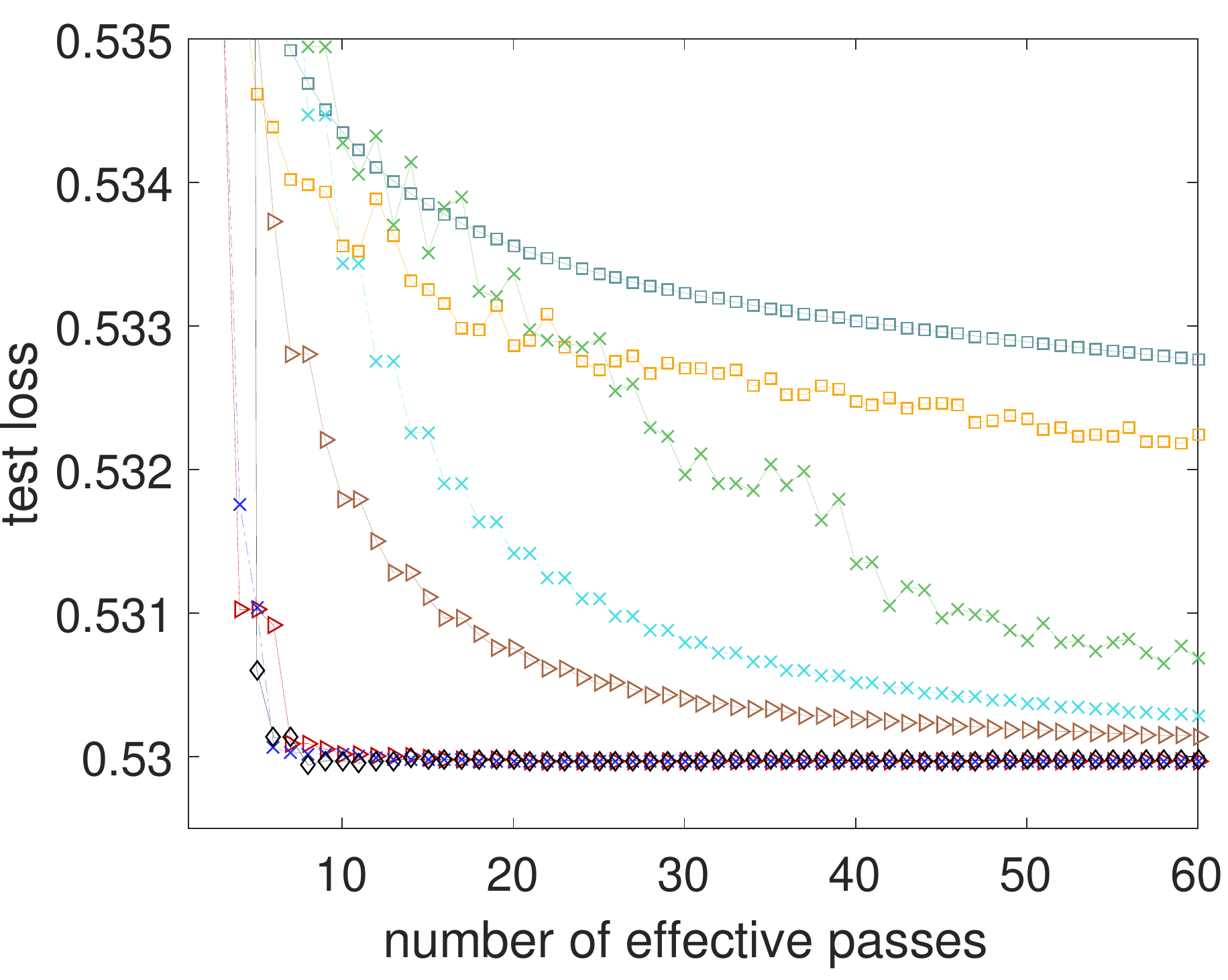}  
	}    
	\subfigure[mnist-testing]{ 
		\label{fig:7} 
		\includegraphics[width=3.4cm,height=3cm]{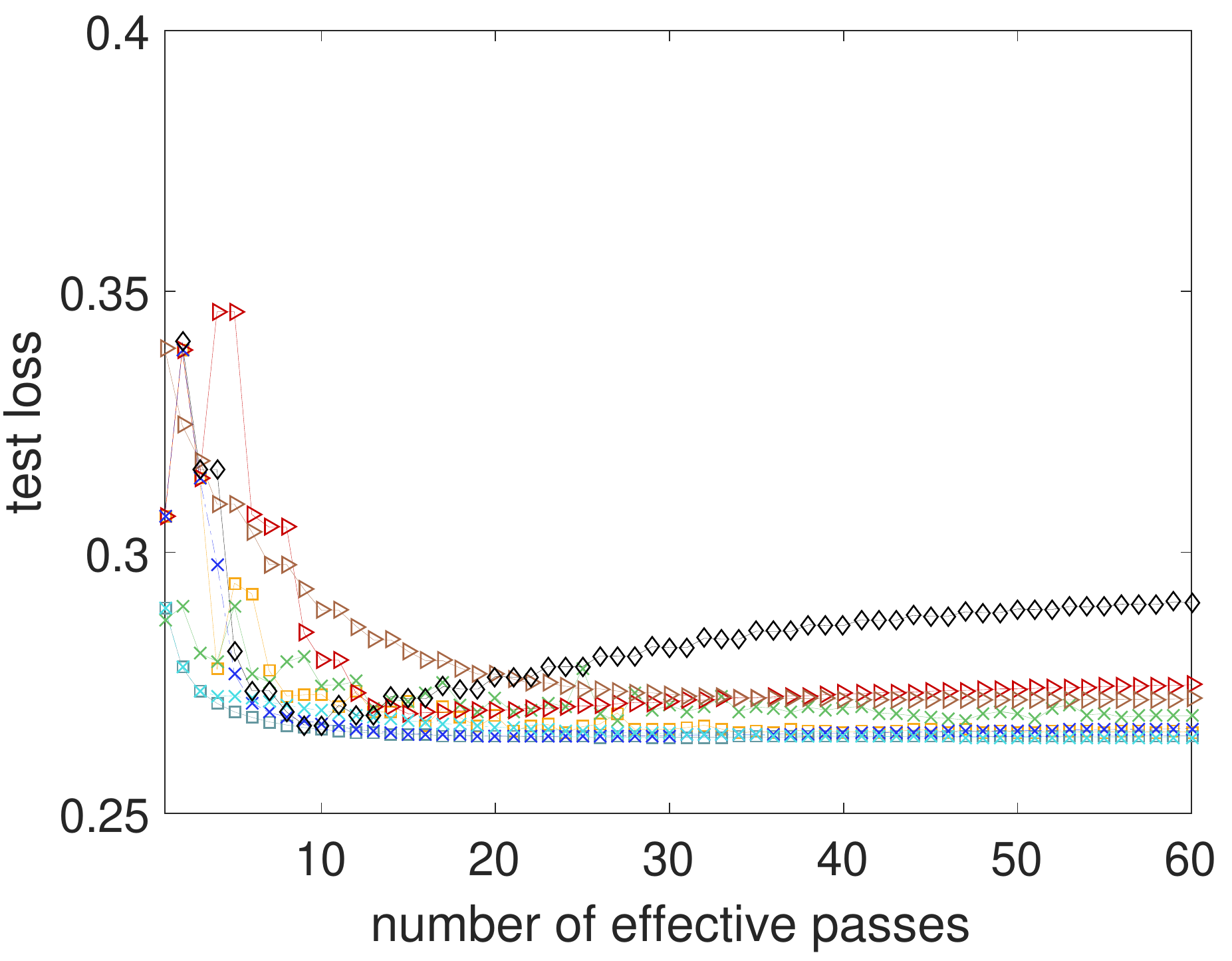}}  
	\subfigure[dna-testing]{ 
		\label{fig:8} 
		\includegraphics[width=3.4cm,height=3cm]{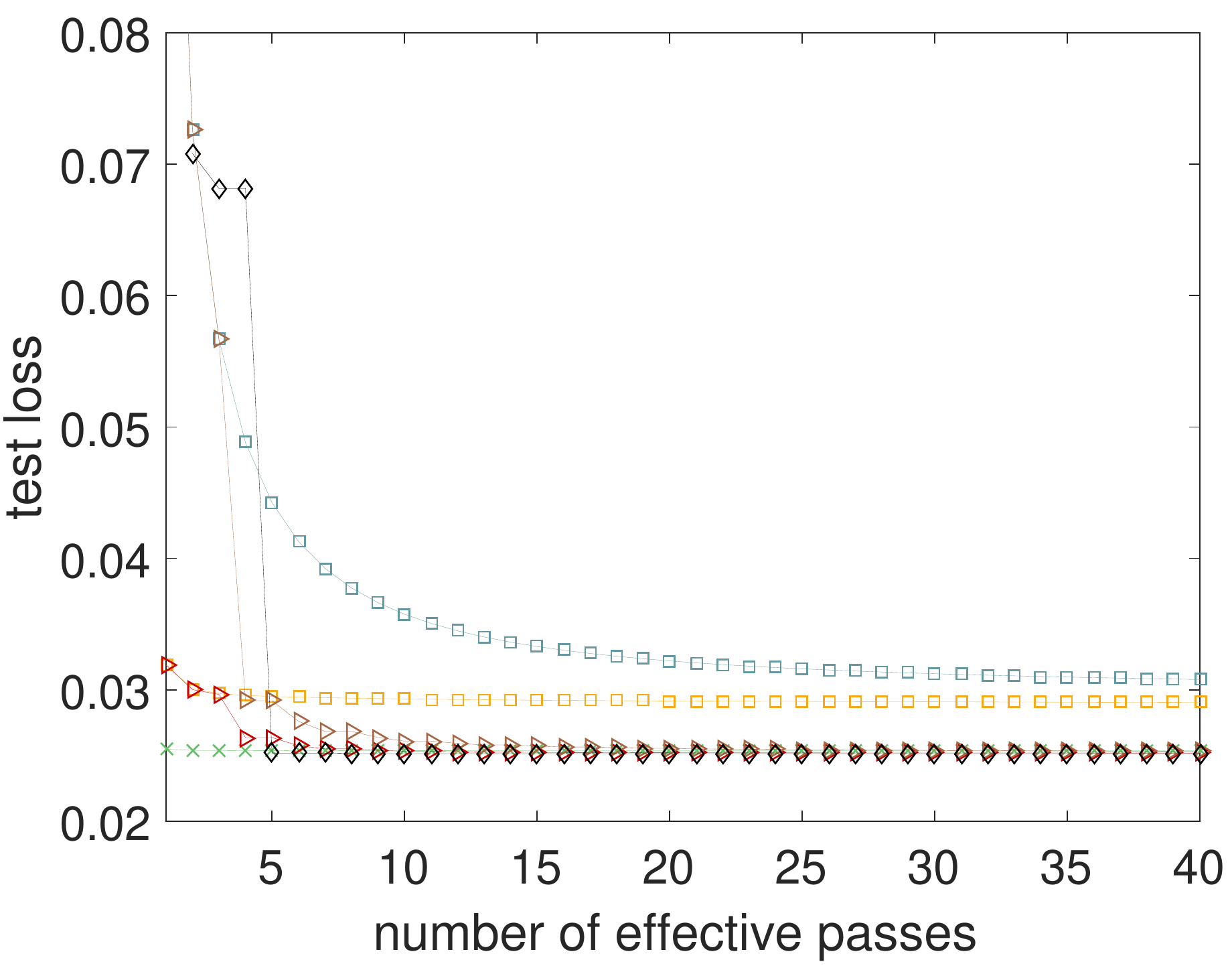}}  	\vspace{-0.1in}
	\subfigure{
		\label{fig:9} 
		\includegraphics[width=0.9\linewidth]{pic/tag.pdf}}  
	\vspace{-0.1in}     
	\caption{Experimental results of solving the  Graph-Guided Fused Lasso problem (Eq.~\eqref{lass2})~(Top: objective gap; Bottom: testing loss). The step size is tuned to be the best for each algorithm. The computation time  has included the cost of calculating full gradients for SVRG based methods. SVRG-ADMM and SAG-ADMM are initialized by running STOC-ADMM for $\frac{3n}{b}$ iterations. ``-ERG'' represents the ergodic results for the corresponding algorithms. }
	\label{graph-guided fused lasso} 
\end{figure}

SAG-ADMM~\cite{SAG-ADMM} and SVRG-ADMM~\cite{SVRG-ADMM} accelerate  SADMM to ergodic  $O(1/K)$.  In Theorem $9$ of~\cite{LADM-NE},  the authors generate a class of functions showing that the original ADMM has a  tight   non-ergodic $O(1/\sqrt{K})$ convergence rate. When $n=1$,  SAG-ADMM  and SVRG-ADMM are the same as batch ADMM,  so their convergence rates  are no better than $O(1/\sqrt{K})$.  This shows that our algorithm has a faster convergence rate than VR based SADMM methods in non-ergodic sense.  One may deem  that judging convergence in ergodic or non-ergodic is unimportant in practice. However, in experiments, our algorithm is much faster than  VR based SADMM methods. Actually, though VR based SADMM methods have provably  faster rates than STOC-ADMM, the improvement in practice is  evident only iterates are close to the convergence point, rather than at the early stage. Both \citet{SAG-ADMM} and \citet{SVRG-ADMM} show that SAG-ADMM and SVRG-ADMM are sensitive to the initial points. We also find that if the step sizes are set based on the their theoretical guidances, sometimes they are even slower than STOC-ADMM~(see Fig.~\ref{original lasso}) as the early stage lasts longer when the step size is small. Our algorithm is faster than the two algorithms whenever the step sizes are set based on the theoretical guidances or are tuned to achieve the fastest speed~(see Fig.~\ref{graph-guided fused lasso}).  This demonstrates that Nesterov's extrapolation has truly accelerated the speed and the integration of   extrapolation and VR techniques is harmonious and complementary.
\subsection{The Growth of Penalty Factor $\frac{\beta}{\theta_1^s}$}
One can find that the penalty factor $\frac{\beta}{\theta_1^s}$ increases linearly with the iteration. This  might make our algorithm impractical because after dozens of epoches, the large value of  penalty factor might slow down the decrement of  function value. However, in experiments, we have not found any bad influence. There may be two reasons 1. In our algorithm,   $\theta_{1,s}$  decreases after each epoch~($m$ iterations),  which is much slower than LADM-NE~\cite{LADM-NE}. For most stochastic problems, algorithms converge in less than $100$ epoches, thus $\theta_{1,s}$ is only  $200$ times of $\theta_{1,0}$. The growth of penalty factor works as  a continuation technique~\cite{Zuo2011A}, which may help to decrease the function value. 2. From Theorem~\ref{convergence1}, our algorithm converges in $O(1/S)$ whenever $\theta_{1,s}$ is large. So from the theoretical viewpoint, a large $\theta_{1,s}$ cannot slow down our algorithm. We find that OPT-ADMM~\cite{OPT-SADMM} also needs to decrease the step size with the iteration. However, its step size decreasing rate is $O(k^{\frac{3}{2}})$ and is faster than ours. 

\section{Experiments}
We conduct experiments to show the effectiveness of our method\footnote{We will put our code online once our paper is accepted.}. We compare our method with the following the-state-of-the-art SADMM algorithms: (1) STOC-ADMM~\cite{STOC-ADMM}, (2) SVRG-ADMM~\cite{SVRG-ADMM}, (3) OPT-SADMM~\cite{OPT-SADMM}, (4) SAG-ADMM~\cite{SAG-ADMM}. We ignore the comparison with SDCA-ADMM~\cite{SDCA-ADMM} since there is no analysis for it on general convex problems and it is  also not faster than SVRG-ADMM~\cite{SVRG-ADMM}. Experiments are performed on  Intel(R) CPU i7-4770 @ 3.40GHz machine with $16$ GB memory. 

\subsection{Lasso Problems}
We perform experiments to solve two typical Lasso problems. The first is the original Lasso problem: 
\begin{eqnarray}\label{lasso}
\min_{\x}\mu\| \x \|_1 +\frac{1}{n}\sum_{i=1}^n\|h_i-\x^Ta_i \|^2,
\end{eqnarray}
where $\mathbf{h}_i$ and $\mathbf{a}_i$ are the tag and the data vector, respectively.  The second is Graph-Guided Fused Lasso model:
\begin{eqnarray}\label{lass2}
\min_{\x}\mu\| \A\x \|_1 +\frac{1}{n}\sum_{i=1}^n l_i(\x),
\end{eqnarray}
where $l_i(\x)$ is the logistic loss on sample $i$, and $\A = [\G; \I]$ is a matrix encoding the feature sparsity pattern.  $\G$ is the sparsity pattern of the graph obtained by sparse inverse covariance estimation~\cite{friedman2008sparse}. The experiments are performed on four benchmark data sets: a9a, covertype, mnist and dna\footnote{a9a, covertype and dna are from:  \url{http://www.csie.ntu.edu.tw/~cjlin/libsvmtools/datasets/}, and mnist is from: \url{http://yann.lecun.com/exdb/mnist/}.}.  The details of the dataset and the mini-batch size that we use in all SADMM are shown in Table~\ref{dataset}.  And like~\cite{SAG-ADMM} and ~\cite{SVRG-ADMM}, we fix $\mu$ to be $10^{-5}$ and  report the performance based on $(\x_t, \A\x_t)$ to satisfy the constraints of ADMM. Results are averaged over five repetitions.  And we  set $m=\frac{2n}{b}$ for all the algorithms. To solve the problems by ADMM methods, we introduce an variable $\y=\x$ or $\y=\A\x$. The update for $\x$ can be written as:  $\x^{k+1} =  \x^{k}-\gamma(\tna f_2(\x_2)+\beta \A^T(\A\x-\y) +\A^T\olam )$, where $\gamma$ is the step size, which depends on the penalty factor $\beta$ and the Lipschitz constant    $L_2$.  For the original Lasso (Eq~\eqref{lasso}), $L_2$ has a closed-form solution, namely, we set $L_2=\max_i \| a_i \|^2=1$\footnote{We normalize the Frobenius norm of each feature to 1.}. So  in this task, the step sizes are set through theoretical guidances  for each algorithm. For Graph-Guided Fused Lasso (Eq~\eqref{lass2}), we regard $L_2$ as a tunable parameter and tune  $L_2$ to obtain the best step size  for each algorithm, which is similar to \cite{SVRG-ADMM} and \cite{SAG-ADMM}.  Except ACC-SADMM, we use the continuation technique~\cite{Zuo2011A} to accelerate  algorithms. We set $\beta_s = \min(10,\rho^s \beta_0)$.   Since SAG-ADMM and SVRG-ADMM are sensitive to initial points,  like \cite{SVRG-ADMM}, they are initialized by running STOC-ADMM for $\frac{3n}{b}$ iterations. SAG-ADMM is performed on the first three datasets due to its large memory requirement.  More details about parameter  setting can be found in  Supplementary Materials.

The experimental results of  original Lasso  and Graph-Guided Fused Lasso are shown in Fig.~\ref{original lasso} and Fig.~\ref{graph-guided fused lasso}, respectively. From the results, we can find that  SVRG-ADMM performs much better than STOC-ADMM when the step size is large while the improvement is not large when the step size is small as it has the cost of calculating the full gradient. SAG-ADMM  encounters a similar situation as $\x$ is not updated  on the latest information. OPT-ADMM performs well on the small step size. However, when the step size is large, the noise of the gradients limits the effectiveness of the extrapolation. Our algorithm is  faster than other SADMM on both  problems.  More experimental results where we set a larger fixed step size  and the memory costs of all algorithms are shown in Supplementary Materials.
\begin{figure}[tb]
	\centering 
	\hspace{-0.15in} 
	\subfigure[objective gap vs. iteration ]{ 
		{\label{fig:1}} 
		\includegraphics[width=0.4\linewidth]{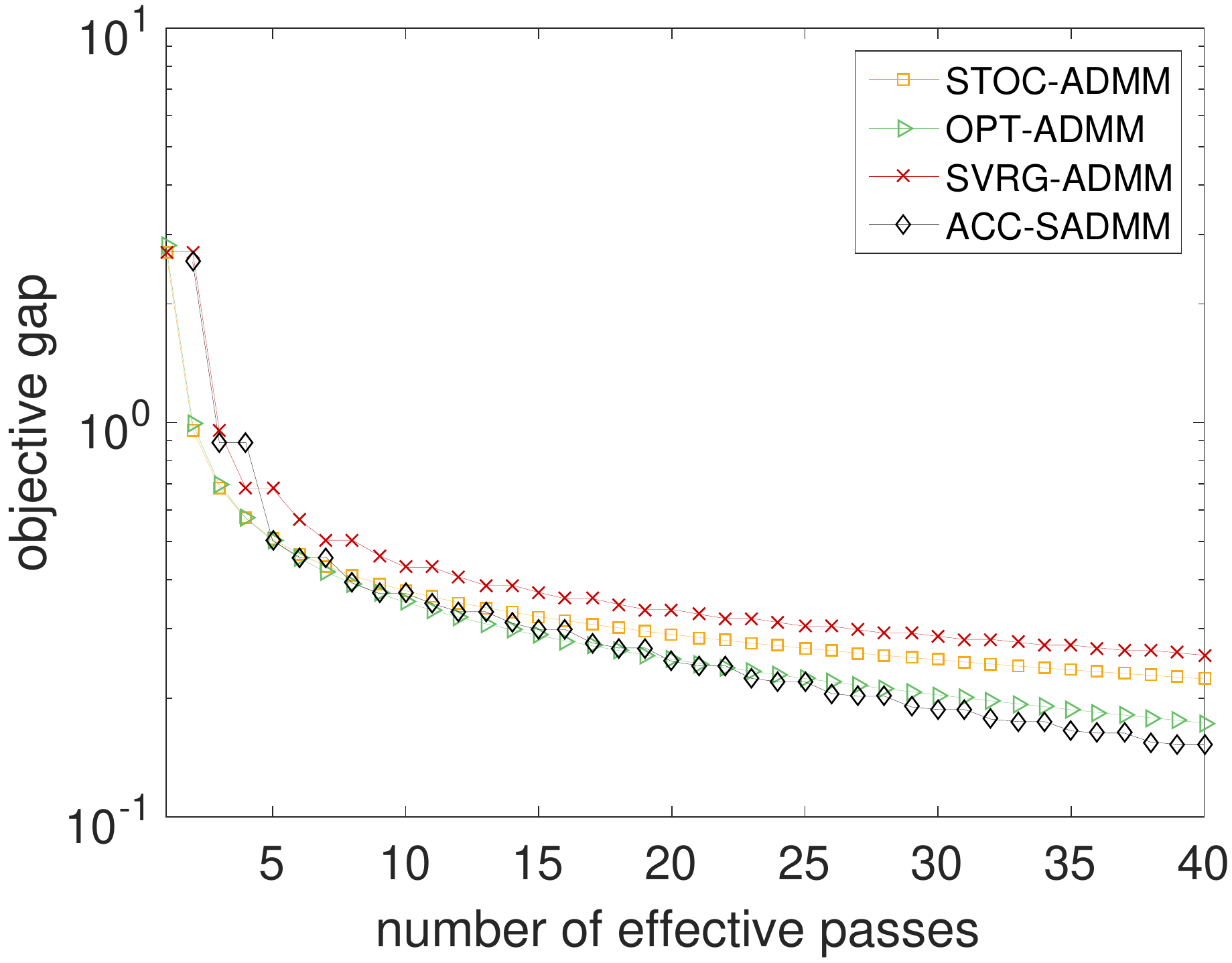}} 
	\subfigure[test error vs. iteration]{ 
		\label{fig:2} 
		\includegraphics[width=0.4\linewidth]{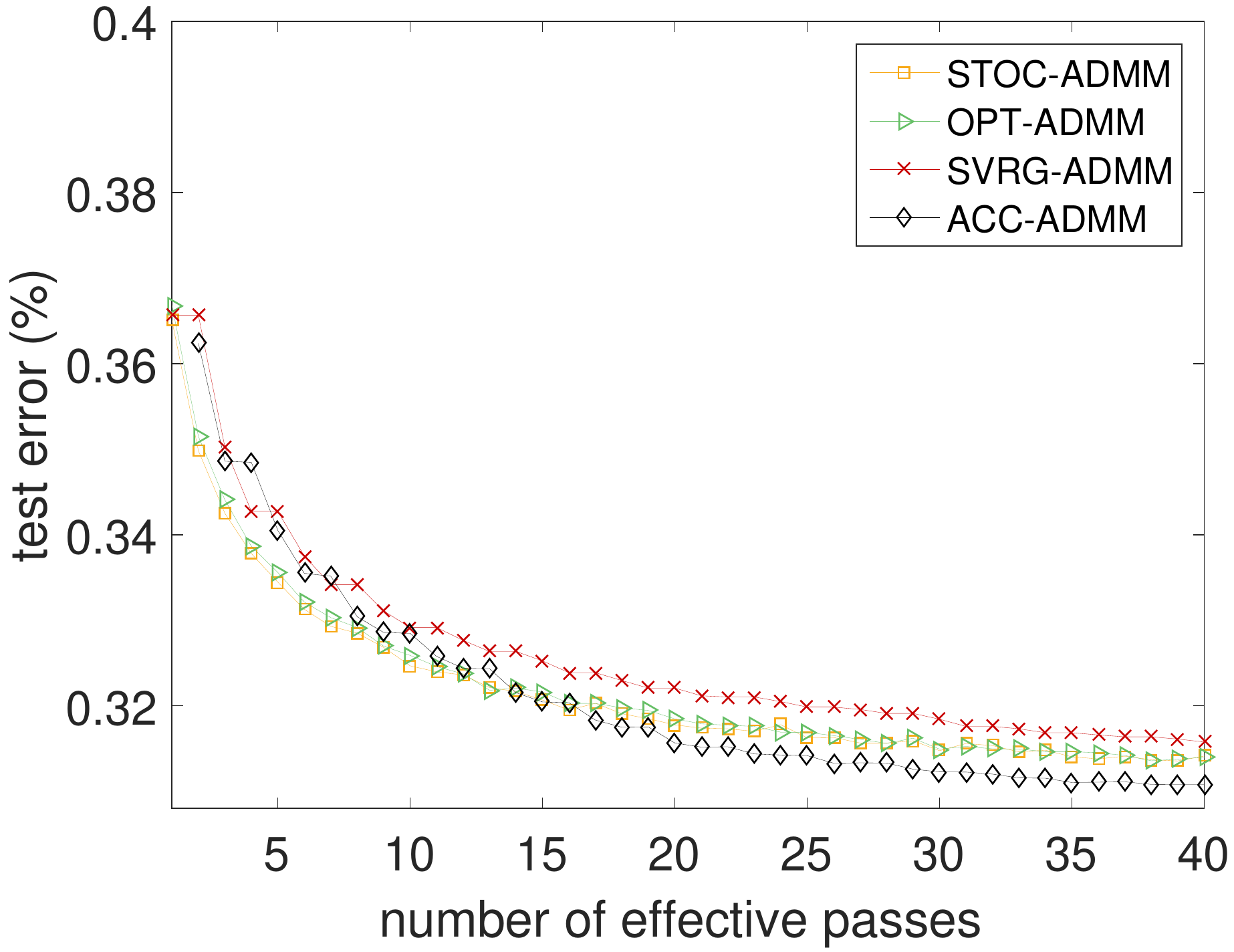}}    
	\caption{The experimental result of Multitask Learning. } 
	\label{vgg1}
\end{figure}

\begin{table}[t]
	\small
	\centering
	\caption{Details of datasets. (Dim., Cla, and mini., are short for dimensionality, class, and mini-batch, respectively. Las. is short for the Lasso problem. Mul. is short for Multitask Learning.).}
	\begin{tabular} {|c |c |c|c|c|c|c|}\hline
		Pro.&Dataset&\# training&\# testing& Dim. $\times$ Cla. &\# mini.\\\hline\hline
		$\!$\multirow{3}*{\centering Las.}$\!\!$&a9a&$72,876$&$72,875$&$74\times2$&\multirow{3}*{\centering $100$}\\\cline{2-5}
		&covertype&$290,506$&$290,506$&$54\times 2$&\\\cline{2-5}
		&mnist&$60,000$&$10,000$&$784\times 10$&\\\cline{2-6}
		&dna&$2,400,000$&$600,000$&$800\times2$&$500$\\\hline
		Mul.&ImageNet&$1,281,167$&$50,000$&$4,096\times 1,000$&$2,000$\\\hline
	\end{tabular}
	\label{dataset}
\end{table}

\subsection{Multitask Learning}
We perform experiments on multitask learning~\cite{multitask}. A similar experiment is also conducted by \cite{SVRG-ADMM}. The experiment is performed on a 1000-class ImageNet dataset~\cite{russakovsky2015imagenet} (see Table~\ref{dataset}). The features are generated from the last fully connected layer of the convolutional VGG-16 net~\cite{simonyan2014very}. More detailed descriptions on the problem are shown in Supplementary Materials.

Fig.~\ref{vgg1} shows the objective gap and test error against iteration. Our method is also faster than other SADMM.

\section{Conclusion}
We  propose  ACC-SADMM  for the general convex finite-sum problems.  ACC-SADMM  integrates   Nesterov's extrapolation and VR techniques and achieves a non-ergodic $O(1/K)$ convergence rate.  We do experiments to demonstrate that our  algorithm is  faster than other SADMM.

\bibliographystyle{custom}
\setlength{\bibsep}{3pt}
\bibliography{ACC_SADMM}


\end{document}